\documentclass[11pt]{article}
\usepackage[colorlinks=true,pdfborder=001,linkcolor=blue,citecolor=blue]{hyperref}
\usepackage{pifont}
\usepackage{url}
\usepackage{stackengine}
\usepackage{nccmath}
\usepackage{mathtools}

\usepackage{mathrsfs}
\usepackage[title]{appendix}
\usepackage{adjustbox}
\usepackage{subcaption}
\usepackage{pgfplots}
\usepackage{tikz}
\usepackage{tikz,fullpage}
\usepackage{rotating}
\usepackage{tabu}
\usepackage{hhline}
\usepackage{multirow}
\usepackage{float}
\restylefloat{table}
\usepackage{adjustbox}
\usepackage[round]{natbib}
\usepackage{booktabs}
\usepackage{array}
\usepackage{tabularx}
\usepackage[utf8]{inputenc}
\usepackage{amssymb,amsmath}
\usepackage{amsfonts}
\usepackage{latexsym}
\usepackage[english]{babel}
\usepackage[T1]{fontenc}
\usepackage{eepic}
\usepackage{epic}
\usepackage{epsfig}
\usepackage{graphics,color}
\usepackage{verbatim}
\usepackage{lscape}
\usepackage{amsthm}
\usepackage{amsmath}
\usepackage{placeins}
\usepackage{float}
\usepackage{caption}
\usepackage{color}
\usepackage{setspace}
\usepackage[ruled,linesnumbered]{algorithm2e}
\usepackage{frcursive}

\usepackage{mathrsfs}
\usepackage{soul}
\usepackage{tikz}
\usepackage{tikz,fullpage}
\usetikzlibrary{arrows,petri,topaths,automata}
\usepackage[short,nocomma]{optidef}
\usepackage{enumitem}  

\usepackage[colorinlistoftodos]{todonotes}

\makeatletter

\renewcommand*\thesection{\arabic{section}}

\newtheorem{prop}{Proposition}

\newtheorem{exam}{Example}

\newif\ifbold
\boldfalse
\boldtrue
\newcommand{\bbf}{\ifbold\bgroup\bf\fi}
\newcommand{\ebf}{\ifbold\egroup\fi}

\renewcommand{\textbf}[1]{\begingroup\bfseries\mathversion{bold}#1\endgroup}
\makeatletter    
\renewcommand{\section}{\@startsection {section}{1}{\z@}%
             {-2ex \@plus -1ex \@minus -.2ex}%
             {1ex \@plus.2ex}%
             {\normalfont\Large\rmfamily\bfseries}}
\renewcommand{\subsection}{\@startsection{subsection}{2}{\z@}%
             {-1.25ex\@plus -1ex \@minus -.2ex}%
             {.75ex \@plus .2ex}%
             {\normalfont\large\rmfamily\bfseries}}
\def\@listI{\leftmargin\leftmargini       
            \parsep .25ex \@plus .1ex     
            \topsep .25ex \@plus .1ex     
            \itemsep \parsep}
\let\@listi\@listI
\@listi
\makeatother
\makeatletter

\@addtoreset{equation}{section}
\makeatother
\makeatletter

\@addtoreset{table}{section}
\makeatother

\usepackage[margin=1in]{geometry}
\usepackage{graphicx}
\graphicspath{ {./Figures/} }
\definecolor{purple}{rgb}{0.4,0.2,1}
\usepackage{titling}

\title{
\LARGE\bf A Rolling-Space Branch-and-Price Algorithm for the Multi-Compartment Vehicle Routing Problem with Multiple Time Windows
\vspace{1ex}
}

\author{\large El Mehdi Er Raqabi,$^{1,2,*}$ Kevin Dalmeijer,$^{1}$ Pascal Van Hentenryck$^{1}$ \\
\footnotesize$^1$\emph{H. Milton Stewart School of Industrial and Systems Engineering, Georgia Institute of Technology, Atlanta, USA}\\
\footnotesize$^2$\emph{Department of Operations and Decision Systems, Université Laval, Québec, Canada}\\
\footnotesize$^*$\emph{Corresponding Author: eraqabi6@gatech.edu}\\
}

\date{}

\begin{document}
\maketitle

\vspace{2cm}
\begin{abstract}
\vspace{0.5cm}

This paper investigates the multi-compartment vehicle routing problem with multiple time windows (MCVRPMTW), an extension of the classical vehicle routing problem with time windows that considers vehicles equipped with multiple compartments and customers requiring service across several delivery time windows. The problem incorporates three key compartment-related features: (i) compartment flexibility in the number of compartments, (ii) item-to-compartment compatibility, and (iii) item-to-item compatibility. The problem also accommodates practical operational requirements such as driver breaks. To solve the MCVRPMTW, we develop an exact branch-and-price (B\&P) algorithm in which the pricing problem is solved using a labeling algorithm. Several acceleration strategies are introduced to limit symmetry during label extensions, improve the stability of dual solutions in column generation, and enhance the branching process. To handle large-scale instances, we propose a rolling-space B\&P algorithm that integrates clustering techniques into the solution framework. Extensive computational experiments on instances inspired by a real-world industrial application demonstrate the effectiveness of the proposed approach and provide useful managerial insights for practical implementation.
 
\vspace{0.2cm}
\par
\vspace{0.2cm}

\end{abstract}    

\setlength{\parindent}{1em}
\setlength{\parskip}{0.5em}
\doublespacing

\newpage

\section{Introduction}\label{section:1}

In today’s rapidly evolving logistics landscape, delivery operations must accommodate increasingly diverse and demanding customer requirements. Many organizations face the challenge of distributing multiple incompatible products that cannot be mixed (e.g., temperature-sensitive product types such as frozen, chilled, and dry goods in grocery distribution, or different fuel types in energy logistics). This led companies to shift gradually toward multi-compartment vehicles rather than traditional single-compartment vehicles. Consequently, the multi-compartment vehicle routing problem has gained both practical and methodological relevance (see e.g., \cite{ostermeier2021multi} and \cite{banipetrol2}). Multi-compartment structures enable fleet size reduction, lower emissions, efficient capacity management, and support for mixed-load operations, which are essential for achieving sustainable and cost-efficient logistics.

Concurrently, the rapid growth of e-commerce and the surge in parcel-based home deliveries have intensified the need for precise, customer-oriented last-mile logistics. In response, many companies now offer customers a selection of narrow, predefined delivery time windows, allowing greater control over delivery timing. While this approach introduces significant constraints on operational planning and routing flexibility, it is widely adopted as a strategic necessity to enhance customer satisfaction, maintain competitiveness, secure shipment volumes from major e-commerce partners, and minimize the costs associated with failed or repeated deliveries. As a result, the vehicle routing problem with multiple time windows has also attracted increasing research and practical attention (see e.g., \cite{schaap2022large}).

Against this background and motivated by a large-scale real-world case, this paper investigates the multi-compartment vehicle routing problem with multiple time windows (MCVRPMTW). The MCVRPMTW extends the traditional VRP with time windows to vehicles with partitioned loading areas (i.e., independent spaces referred to as compartments) and customers who may have multiple delivery time windows. While each of these components (multi-compartment structures and multiple time windows) independently introduces substantial operational and computational challenges, their joint consideration reflects the real-world complexity faced by modern logistics organizations. Addressing both simultaneously enables the development of solution approaches that support efficient, service-oriented, and resource-conscious delivery operations aligned with current industry demands for multi-product transport and precisely timed customer service.

The MCVRPMTW considered in this study is motivated by a large-scale problem provided by our industrial partner. This real-world case involves loading constraints on compartmentalized vehicles, geographically dispersed urban customers, multiple delivery options across a weekly planning horizon, and additional operational constraints such as driver regulations. Beyond addressing this specific industrial challenge, the objective of this paper is to propose a general methodological framework for solving the MCVRPMTW. The framework is designed to be both practically applicable and methodologically extensible, advancing the frontier of VRP-related research. 

The MCVRPMTW inherits the computational complexity of several NP-hard components. The classical VRP with time windows (VRPTW) is itself NP-hard \citep{solomon1987algorithms}, and the addition of multiple compartments introduces a further combinatorial layer arising from loading restrictions and compatibility constraints \citep{ostermeier2018loading}. Each vehicle configuration may yield an exponentially large number of feasible loading patterns, while time-window synchronization tightly couples routing and scheduling decisions. Consequently, the MCVRPMTW integrates spatial, temporal, and loading complexities, rendering it substantially more challenging than traditional single-compartment VRP variants.

To the best of our knowledge, this work is the first to formally introduce and study the MCVRPMTW, a problem motivated by practice and offering methodological novelty. The contributions are fourfold:
\begin{itemize}
    \item The paper defines a new routing framework that jointly incorporates multiple vehicle compartments and multiple customer time windows. The compartment-related aspects explicitly captured include compartment flexibility, item–compartment compatibility, and item–item compatibility \citep{ostermeier2021multi, cherkesly2022pickup}. The formulation further integrates operational elements observed in practice, such as driver break requirements.
    \item The paper models the MCVRPMTW as a set-partitioning problem and designs an exact branch-and-price (B\&P) algorithm in which the pricing subproblem is solved via a labeling procedure. The paper proposes two theoretically valid dominance rules, \emph{same-day} and \emph{inter-day} dominance rules, that maintain correctness while substantially reducing the number of labels. Several acceleration strategies are incorporated to mitigate symmetry in label extensions, stabilize dual solutions during column generation, strengthen branching decisions, and efficiently manage the breaks.
    \item To address large-scale instances, the paper proposes a rolling-space B\&P algorithm that integrates clustering techniques into the solution framework.
    \item The paper conducts extensive computational experiments on a large dataset derived from an industrial application, demonstrating the effectiveness and scalability of the proposed approach and offering actionable managerial insights for practical implementation.
\end{itemize}

\noindent
The remainder of this paper is organized as follows. 
Section~\ref{section:2} reviews related literature. Section~\ref{section:3} presents the MCVRPMTW formulation. Section~\ref{section:4} details the B\&P algorithm. Computational results are reported in Section~\ref{section:5}, and Section~\ref{section:6} concludes the paper.

\section{Literature Review}\label{section:2}

The MCVRPMTW falls within the broader class of VRPs. It generalizes both the VRP with time windows and the multi-compartment VRP (see recent surveys by \cite{ostermeier2021multi}, \cite{mor2022vehicle}, and \cite{archetti2025beyond}). In what follows, we first review the literature on multi-compartment VRPs and then summarize contributions related to VRPs with multiple time windows.

\subsection{Multi-Compartment VRP}

Multi-compartment vehicle routing problems (VRPMCs) extend classical VRPs by allowing vehicles to carry goods in several distinct compartments. Such models arise in several real-world contexts, including petrol distribution, maritime logistics, and waste collection (see \cite{coelho2015classification} and \cite{ostermeier2021multi} for detailed surveys). Prior research has mainly examined three compartment-related dimensions: (i) compartment size flexibility, (ii) item-to-compartment compatibility, and (iii) item-to-item compatibility \citep{ostermeier2021multi}.

Most contributions in the VRPMC literature consider only one or two compartment-related features at a time \citep{ostermeier2018loading}. Many studies assume fixed compartment capacities (e.g., \cite{coelho2015classification} and \cite{zbib2020commodity}), although some works allow compartment sizes to vary (e.g., \cite{henke2015multi}, \cite{hubner2019multi}, and \cite{hessler2023partial}). Item–compartment assignments may also be either restricted, where each product type can be loaded only in a prescribed subset of compartments (e.g., \cite{ostermeier2018loading} and \cite{martins2019product}), or flexible, allowing items to use any compartment (e.g., \cite{lahyani2015multi} and \cite{christiansen2017operational}). Another modeling choice concerns whether compartments can serve multiple customers: most studies allow shared compartments (e.g., \cite{kiilerich2018new} and \cite{yahyaoui2020two}), whereas a smaller number examine settings in which each compartment is dedicated exclusively to a single customer (e.g., \cite{jetlund2004improving} and \cite{hsu2020solving}).

Several exact B\&P algorithms have been proposed for different variants of the VRPMC, distinguished by elements such as the number of compartments, product–compartment assignment rules, and computational scale. \cite{avella2004solving} study a variant with fixed compartment capacities, flexible product-to-compartment assignments, and unshared compartments. Their route-enumeration B\&P algorithm efficiently solves real-world instances with up to 60 customers, using vehicles with 7–9 compartments in only a few seconds. \cite{mirzaei2019branch} consider a model with fixed capacities, fixed product assignments, and shared compartments, distinguishing between single- and multiple-visit customers. Their method handles instances with up to 100 customers and 4 compartments, and proves optimality for instances of up to 50 customers. \cite{hessler2021exact} examines a variant with flexible compartment capacities and constrained product-category-to-compartment assignments. Using a branch-and-price-and-cut algorithm, they solve instances with up to 50 customers and 2–9 compartments within two hours. More recently, \cite{hessler2023partial} propose two labeling-based B\&P algorithms that handle similar characteristics but focus on computational efficiency through dominance rules. Their method handles instances with up to 100 customers and 4 compartments, solving them optimally within one hour. \cite{bani2023solving} address a new variant of the petrol replenishment problem, which is a rich real-world multi-depot multi-period problem, using B\&P combined with several acceleration techniques that exploit the polyhedral properties of the problem. They solve many real-world network (four depots, five types of petroleum products, four main groups of clients, and a heterogeneous fleet of highly compartmented tank trucks) instances to near-optimality. In a follow-up extension, \cite{bani2024combining} also incorporate the inventory routing problem. Driven by a decoupling intuition, the authors develop an exact two-phase solution approach that combines Benders decomposition and column generation. In the first phase, they solve the relaxed (integrality) Benders subproblems using column generation until the inventory levels stabilize. In the second phase, they solve the Benders subproblems using column generation embedded in a branch-and-bound (B\&B) framework. They reach near-optimal solutions in all these instances and note that acceleration strategies significantly boost the performance of the two-phase method. Lastly, \cite{aerts2024unified} present a unified B\&P algorithm for multi-compartment pickup and delivery problems. The algorithm is evaluated on instances featuring up to 75 requests, 3 compartments, and 24 item categories, successfully solving 90\% of the instances within 3,600 seconds.

\subsection{VRP with Multiple Time Windows}

Compared to the VRPTW and the MCVRP, the VRPMTW has received considerably less attention in the literature. In the MCVRPMTW, each customer may be associated with multiple alternative time windows, representing different feasible delivery periods. However, each customer must be served exactly once, within only one of these time windows, distinguishing the problem from split-delivery variants. \cite{favaretto2007ant} provided one of the first formal definitions of the VRPMTW and examined two objective functions: minimizing total route duration and minimizing total travel distance. They proposed an ant colony optimization algorithm and generated new benchmark instances derived from the VRP dataset of \cite{fisher1994optimal}. Building on this work, \cite{belhaiza2014hybrid} developed a hybrid variable neighborhood tabu search that achieved improved solution quality and introduced additional benchmark instances inspired by customer patterns from \cite{solomon1987algorithms}. \cite{amorim2014rich} applied an adaptive large neighborhood search (ALNS) to a real-world case study involving perishable food distribution, considering a site-dependent VRP with a heterogeneous fleet.

Further advancements were achieved by \cite{belhaiza2017new}, who improved upon previous results through a genetic variable neighborhood search. In subsequent work, \cite{belhaiza2019three} reported numerous new best-known solutions (BKSs) for the route duration minimization objective and demonstrated the applicability of their approach in a real-world Canadian furniture delivery case. In parallel, \cite{larsen2019fast} proposed an ALNS algorithm for the VRPMTW that incorporates efficient cost update mechanisms, conceptually related to those employed in the present study. However, their work does not provide a rigorous analysis of the correctness and computational complexity of the underlying dynamic programming components and offers a limited comparative discussion, focusing primarily on \cite{belhaiza2014hybrid}.

Additional contributions include \cite{hoogeboom2020efficient}, who developed an adaptive variable neighborhood search (AVNS) that incorporates an efficient algorithm for evaluating route durations on fixed customer sequences. \cite{ferreira2018variable} also proposed a variable neighborhood search combined with a column generation approach, solving VRPMTW instances with up to 17 customers to optimality, as later cited by \cite{bogue2022column}. More recently, \cite{schaap2022large} introduced a large neighborhood search–based metaheuristic for the VRPMTW. Their algorithm was tested on the benchmark set of \cite{belhaiza2014hybrid}, yielding new best-known solutions for 9 out of 48 instances for distance minimization and 13 out of 48 instances for duration minimization.

Constraints involving multiple time windows have been widely explored in related problem settings, including the traveling salesperson problem, truck driver scheduling, and the team orienteering problem. \cite{xu2003solving} examined a pickup-and-delivery problem with time windows and addressed it using a heuristic column-generation approach. \cite{tricoire2010heuristics} studied the team orienteering problem with multiple time windows, proposing a variable-neighborhood search that embeds an exact time-window-allocation subproblem. \cite{paulsen2015heuristic} investigated the traveling salesperson problem with multiple time windows (TSPMTW), evaluating a heuristic dynamic-programming method and a genetic algorithm to optimize customer sequences and departure times jointly. \cite{baltz2015exact} extended the TSPMTW to incorporate hotel selection and overnight stays, introducing a mixed-integer programming formulation combined with a randomized heuristic.

\section{Mathematical Formulation}\label{section:3}

The MCVRPMTW is modeled on a directed graph $\mathcal{G}=(\mathcal{V},\mathcal{A})$, where $\mathcal{V}=\mathcal{N} \cup \{\sigma,\delta\}$ represents the set of vertices and $\mathcal{A}$ the set of arcs. The depot is represented by the origin $\sigma$ and destination $\delta$, while $\mathcal{N}$ denotes the set of clients. Each client $n \in \mathcal{N}$ is characterized by a demand $\boldsymbol{Q}_n$, a service duration $\boldsymbol{u}_n$, and a set of $\boldsymbol{\theta}_n$ allowable time windows $\Theta_n = \{[\boldsymbol{s}_n^1,\boldsymbol{e}_n^1], \dots, [\boldsymbol{s}_n^{\theta_n},\boldsymbol{e}_n^{\theta_n}]\}$ for service initiation. Without loss of generality, these time windows are assumed to be disjoint, ordered chronologically, and nonnegative; overlapping windows are merged and sorted for computational purposes, with scheduled arrival times later mapped back to the original windows.

The planning horizon spans multiple days, and each allowable time window is assigned to a specific day $d \in \mathcal{D}$. This multi-day structure reflects the operational reality that customers may be served on different days under distinct availability periods, providing a natural framework for incorporating driver regulations. In particular, for each day $d \in \mathcal{D}$, a driver may not exceed a maximum driving distance $\boldsymbol{Dist}$ and a maximum working time $\boldsymbol{T}$, both of which include travel and service durations. At the end of each daily shift, the driver must take a mandatory rest period of length $\boldsymbol{B}$ before beginning operations on the next day.

Each arc $(i,j) \in \mathcal{A}$ is associated with a travel cost $\boldsymbol{c}_{ij}$ and a travel duration $\boldsymbol{t}_{ij}$, and both are assumed to satisfy the triangle inequality. Each vehicle $k \in \mathcal{K}$ has a total capacity $\boldsymbol{L}_k$ and a set of compartments $\mathcal{M}_k$, where each compartment $m \in \mathcal{M}_k$ has a maximum capacity $\overline{\boldsymbol{L}}_{mk}$ such that $0 \leq \overline{\boldsymbol{L}}_{mk} \leq \boldsymbol{L}_k$. Compatibility between clients and compartments is specified by the binary parameter $\boldsymbol{b}_{nm}^k$, equal to 1 if client $n$ can be assigned to compartment $m$ of vehicle $k$, and 0 otherwise. Similarly, item-to-item compatibility between clients $i$ and $j$ is represented by $\boldsymbol{f}_{ij}$, set to 1 if the items are compatible and 0 otherwise. 

The MCVRPMTW seeks a set of minimum-cost elementary routes from $\sigma$ to $\delta$ that together serve all clients exactly once and satisfy time window constraints, item-to-item compatibility, compartment and vehicle capacities, item-to-compartment compatibility, as well as driver regulations. The following example will illustrate two feasible routes.

\begin{figure}[t!]
    \centering
    \includegraphics[width=0.8\linewidth]{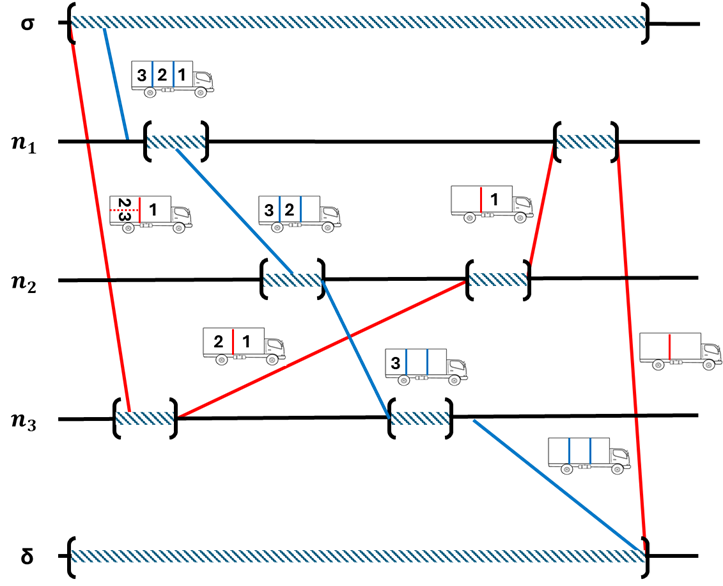}
    \caption{MCVRPMTW Example}
    \label{MCVRPMTW Example}
\end{figure}

\begin{exam}
Consider the MCVRPMTW example in Figure~\ref{MCVRPMTW Example}.
The horizontal layers labeled $n_1$, $n_2$, and $n_3$ represent three customers, each with two possible time windows.
The top and bottom rows correspond to the depot, indicating the start and end of each route.
We assume that the product corresponding to Customer 3 can only be delivered from the back of the vehicle, which is captured through item-to-compartment compatibilities.
Furthermore, the product corresponding to Customer 3 is only compatible with the product corresponding to Customer 2, which is captured through item-to-item compatibilities.\\
The red and blue routes depict two alternative ways to serve the three customers while respecting delivery time and compartmental constraints.
In the blue vehicle with three compartments, Customer 3's order, which can only be delivered from the back, is assigned to the third compartment. In the red vehicle, having only two compartments, the second compartment is shared between Customers 2 and 3.
The two routes differ significantly, yet both respect the multi-compartment and time window constraints.\\
Figure \ref{MCVRPMTW Example} also illuminates the enormous solution space inherent to MCVRPMTW. With enough time windows and vehicle types, there are numerous ways to configure vehicle routes, assign products to compartments, and schedule visits within allowable time windows. Even for a small set of customers, the number of feasible combinations is substantial.
This combinatorial explosion underpins the complexity of the MCVRPMTW and motivates the need for advanced optimization methods such as decomposition techniques.
\end{exam}

Let $\Omega$ denote the set of all feasible routes as described above. Each route $r \in \Omega$ has an associated total cost $\boldsymbol{c}_r$, which includes a fixed vehicle cost plus the sum of travel costs along its constituent arcs. For each client $n \in \mathcal{N}$, let $\boldsymbol{a}_{nr}$ be a binary parameter equal to 1 if route $r$ visits client $n$, and 0 otherwise. The model uses a binary decision variable $y_r$, set to 1 if route $r$ is selected in the solution and 0 otherwise. With these definitions, the MCVRPMTW can be formulated as follows

\begin{align}
\label{Set Partitionning Objective}
\min_{y} \hspace{1mm} & \sum_{r \in \Omega} \boldsymbol{c}_ry_r &&\\
s.t. \hspace{1mm} \label{One Route per Client} & \sum_{r \in \Omega} \boldsymbol{a}_{nr}y_r=1, && \forall n \in \mathcal{N},\\
\label{Binary Condition on Routes} & y_r \in \{0,1\}, && \forall r \in \Omega.
\end{align}

\noindent
The objective function (\ref{Set Partitionning Objective}) minimizes the overall routing costs. Constraints (\ref{One Route per Client}) guarantee that every client is visited exactly once. The conditions on the decision variables are specified by Constraints (\ref{Binary Condition on Routes}).
When travel costs and times obey the triangle inequality, Constraints (\ref{One Route per Client}) can alternatively be expressed as covering constraints, and the binary restrictions (\ref{Binary Condition on Routes}) can be relaxed to nonnegative integer conditions, $y_r \in \mathbb{N}_0$, for all $r \in \Omega$.
Note that the route definition implicitly incorporates time window constraints, item-to-item compatibility, compartment and vehicle capacity limits, item-to-compartment compatibility, and driver regulations.
Appendix \ref{A} describes a compact formulation that explicitly models these aspects.

\section{Branch-and-Price Algorithm}\label{section:4}

Model (\ref{Set Partitionning Objective}) typically has a vast number of variables, corresponding to all possible routes. This makes full enumeration impractical and sometimes intractable. The paper therefore proposes a B\&P approach \citep{barnhart1998branch} for the MCVRPMTW. A B\&P algorithm is a B\&B framework that uses column generation to generate variables and repeatedly solves the linear relaxation of Model (\ref{Set Partitionning Objective}), referred to as the master problem (MP), for an increasing number of variables. The column generation alternates between solving the RMP and a pricing problem to find columns with negative reduced cost that can be added to the RMP. When no improving columns exist, the optimal solution of the relaxed MP is obtained. To ensure integrality, the B\&P algorithm branches on arcs (whether an arc is present in the optimal solution) and performs column generation at each node. The resulting B\&B tree is explored using a best-bound-first search strategy. Other branching rules, such as those based on vehicle counts, can also be employed. For conciseness, see \cite{ropke2009branch} for a detailed overview of this process.

The rest of this section describes the pricing problem of the B\&P algorithm for the MCVRPMTW. Recognizing that classical B\&P methods do not scale well to large instances \citep{bouvier2024solving}, the section then introduces a rolling-space B\&P algorithm that integrates clustering techniques to handle large-scale MCVRPMTW scenarios effectively. Finally, the section discusses the acceleration strategies incorporated into the labeling process.

\subsection{Pricing Problem}\label{Exact B&P}

Let $\pi_n$ denote the dual variable corresponding to Constraint (\ref{One Route per Client}) for each $n \in \mathcal{N}$. The reduced cost of a route $r \in \Omega$ is defined as

\begin{equation}
    \overline{c}_r=c_r-\sum_{n \in \mathcal{N}}\boldsymbol{a}_{nr}\pi_n,
\end{equation}

\noindent i.e., the reduced cost equals the original route cost minus the sum of the dual values of the clients served by the route. The pricing problem involves either identifying a feasible route with negative reduced cost or proving that no such route exists. This problem corresponds to a variant of the elementary shortest path problem with resource constraints (ESPPRC), extended to include multiple time windows, vehicle capacities, multiple compartments with individual maximum capacities, and item-to-compartment as well as item-to-item compatibility constraints. 
The pricing problem is solved by a labeling algorithm \citep{irnich2005shortest}, where partial paths are represented as labels and extended along network arcs. Labels that cannot lead to an optimal solution are discarded using dominance rules.

The proposed labeling algorithm for the MCVRPMTW simultaneously accounts for multiple time windows and the three-compartment features: compartment flexibility, item-to-compartment compatibility, and item-to-item compatibility. The presentation begins by defining the reduced cost and specifying the resources associated with forward labels. It then introduces resource extension functions (REFs) that generate only feasible label extensions with respect to both the time windows and multi-compartment constraints. Finally, two valid dominance rules to prune the set of labels efficiently are presented.

\subsubsection{Reduced Cost}

Let us define $\overline{\pi}_n := \pi_n$ if $n \in \mathcal{N}$, $\overline{\pi}_\sigma := 0$, and $\overline{\pi}_\delta := 0$. Using this definition, the reduced cost of an arc $(i,j) \in \mathcal{A}$ is
\begin{equation}
    \overline{c}_{ij} := c_{ij} - \overline{\pi}_i,
\end{equation}
Accordingly, the reduced cost of a route $r \in \Omega$ is the sum of the reduced costs of its arcs:
\begin{equation}
    \overline{c}_r = \sum_{(i,j) \in \mathcal{A}(r)} \overline{c}_{ij},
\end{equation}
where $\mathcal{A}(r)$ represents the ordered set of arcs forming route $r$.

\subsubsection{Label Resources}

A partial path $R=(\sigma,\dots,i)$, which starts at the depot, is served by vehicle $k$, and ends at client $i$, is represented by a label $L = \bigl(i, t, \Bar{c}, d, \mathcal{S}, (l_m^k)_{m \in \mathcal{M}_k}, (O_m^k)_{m \in \mathcal{M}_k}, g, o\bigr)$, which stores the following information:

\begin{itemize}
    \item $i(L) \in \mathcal{N} \cup \{\sigma,\delta\}$: the latest visited node;
    \item $t(L) \in \mathbb{R}_{+}$: the accumulated time; 
    \item $\Bar{c}(L) \in \mathbb{R}$: the accumulated reduced cost;
    \item $d(L) \in \mathcal{D}$: the day;
    \item $\mathcal{S}(L) \subseteq \mathcal{N}$: the set of visited clients;
    \item $l_m^k(L) \in \mathbb{R}_{+}$: the load in compartment $m \in \mathcal{M}_k$ of vehicle $k$ ($\lvert M_k \rvert$ values per label);
    \item $O_m^k(L) \subseteq \mathcal{N}$: the clients in compartment $m \in \mathcal{M}_k$ of vehicle $k$ ($\lvert M_k \rvert$ sets per label);
    \item $g(L) \in \mathbb{R}_{+}$: accumulated working time for the current day $d \in \mathcal{D}$ with upper bound $\boldsymbol{T}$;
    \item $o(L) \in \mathbb{R}_{+}$: accumulated driving distance for the current day $d \in \mathcal{D}$ with upper bound $\boldsymbol{Dist}$.
\end{itemize}

Starting from the depot, the initial label is defined as $L_0=\left(\sigma,0,0,0,\emptyset,(0)_{m \in \mathcal{M}_k},(\emptyset)_{m \in \mathcal{M}_k},0,0\right)$.

\subsubsection{Resource Extension Function}

The extension of a label $L$ along an arc $(i,j) \in \mathcal{A}$ may have multiple extensions. Each extension is defined by the choice of a time window, a specific compartment within the vehicle, and whether a break is taken. Let $\mathcal{H}(j)$ denote the set of feasible extensions along arc $(i,j)$ that respect the time windows, the multi-compartment, and the driver regulations.
Each extension $h \in \mathcal{H}(j)$ generates a new label $L_h$ according to the Research Extension Functions (REFs) defined below.
Let $\Bar{p} = \min \left\{ p \,\big|\, t(L_h) \leq \boldsymbol{e}_j^p \right\}$ indicate the earliest reachable time window for the extension, $d_{ij}$ be a binary indicator that equals 1 if there is a break $\boldsymbol{B}$ from node $i$ to node $j$ that causes the deliveries to start on the next day from $j$, and 0 otherwise, and $m^h$ be the compartment used by extension $h$. We assume the break $\boldsymbol{B}$ occurs immediately at $i$, and then the vehicle drives to $j$.
The REFs generate the new label $L_h$ as follows:

{\allowdisplaybreaks
\begin{align}
i(L_h) &:=j \label{REF_Node} \\ 
t(L_h) &:=\max\{t(L)+\boldsymbol{t}_{ij}+\boldsymbol{u}_{i}+\boldsymbol{B}d_{ij},\boldsymbol{s}_j^{\Bar{p}}\} \label{REF_Time} \\
\Bar{c}(L_h) &:=\Bar{c}(L)+\overline{c}_{ij} \label{REF_Cost} \\
d(L_h) &:=d(L)+d_{ij} \label{REF_Day} \\
\mathcal{S}(L_h) &:= \mathcal{S}(L_h) \cup \{j\} \label{REF_Client} \\
l_m^k(L_h) &:=\begin{cases}
l_m^k(L)+L_j & \text{if } m=m^h\\
l_m^k(L), & \text{otherwise}
\end{cases} \label{REF_Capacity} \\
O_m^k(L_h) &:=\begin{cases}
O_m^k(L) \cup \{j\} & \text{if } m=m^h\\
O_m^k(L), & \text{otherwise}\\
\end{cases} \label{REF_Visited}\\
g(L_h) &:=\begin{cases}
g(L)+\boldsymbol{t}_{ij}+\boldsymbol{u}_i & \text{if } d_{ij}=0\\
0, & \text{otherwise}
\end{cases} \label{REF_Work_Day} \\
o(L_h) &:=\begin{cases}
o(L)+\boldsymbol{dist}_{ij} & \text{if } d_{ij}=0\\
0, & \text{otherwise}
\end{cases} \label{REF_Drive_Day} 
\end{align}
}

Equations \eqref{REF_Node}–\eqref{REF_Drive_Day} define the REFs that update the label attributes when extending a partial path along an arc. Equation \eqref{REF_Node} updates the current node of the label to the newly visited node $j$. Equation \eqref{REF_Time} updates the accumulated service start time by incorporating travel and service durations while ensuring feasibility with the selected time window. The reduced cost is updated in \eqref{REF_Cost} by adding the cost of arc $(i,j)$. The day counter in \eqref{REF_Day} determines whether the transition occurs within the same day or moves to the next day. Equation \eqref{REF_Client} records the addition of client $j$ to the set of visited clients. The capacity update in Equation \eqref{REF_Capacity} adjusts the load in the compartment $m$ assigned to the extension $h$, while Equation \eqref{REF_Visited} tracks which clients are served by each compartment. Equation \eqref{REF_Work_Day} tracks the daily constraint for working time, while Equation \eqref{REF_Drive_Day} tracks the daily constraint for driving time.

The label $L_h$ is retained if it satisfies the time window and capacity constraints and does not revisit any customer already included in $\mathcal{S}(L_h)$, thereby explicitly prohibiting cycles. That is, a label is retained only if:

\begin{equation}\label{earliest window}
 t(L_h) \leq \boldsymbol{e}_j^{\Bar{p}},
\end{equation}

\begin{equation}\label{limit compartment}
     l_m^k(L_h) \leq \boldsymbol{\overline{L}}_{mk},
\end{equation}

\begin{equation}\label{limit truck}
    \sum_{m\in \mathcal{M}_k} l_m^k(L_h) \leq \boldsymbol{\overline{L}}_k,
\end{equation}

\begin{equation}\label{daily time limit}
 g(L_h) \leq \boldsymbol{T},
\end{equation}

\begin{equation}\label{daily drive limit}
 o(L_h) \leq \boldsymbol{Dist},
\end{equation}

Equation~\eqref{earliest window} implies that the extension selects the earliest possible time window whose end time is later than the updated accumulated time, accounting for service and travel times. Consequently, the updated time for the label is $t(L_h):=\max\{t(L)+\boldsymbol{t}_{ij}+\boldsymbol{u}_{i}+\boldsymbol{B}d_{ij},\boldsymbol{s}_j^{\Bar{p}}\}$ since delaying service to the beginning of its feasible window does not affect route feasibility or cost and is therefore without loss of optimality. The assignment of client $j$ to a compartment $m \in \mathcal{M}_k$ must satisfy both the item-to-compartment compatibility, i.e., the items associated with client $j$ are compatible with compartment $m$, and the item-to-item compatibility, i.e., client $j$ is compatible with all items already loaded in compartment $m$, as well was the capacities in the compartment and the truck (Equations~\eqref{limit compartment} and \eqref{limit truck}). Equations~\eqref{daily time limit} and \eqref{daily drive limit} control the driver regulations regarding the daily constraints.

\subsubsection{Label Dominance}

In the presence of driver breaks and daily resource constraints, label dominance becomes more intricate. To efficiently prune suboptimal labels while ensuring feasibility, two dominance rules are considered. The first, \emph{same-day dominance}, compares labels only within the same day, providing a simple and correct pruning criterion. The second, \emph{inter-day dominance}, allows labels from earlier days to dominate labels from later days if there is enough time to catch up on a break. 

\noindent \textbf{Same-day Dominance.} In this same-day approach, labels are compared only if they belong to the same day. A label $L_1$ dominates another label $L_2$ if

{\allowdisplaybreaks
\begin{subequations}\label{Dominance}
\begin{align}
i(L_1) &= i(L_2), \label{Dominance_a} \\
t(L_1) &\le t(L_2), \label{Dominance_b} \\
\Bar{c}(L_1) &\le \Bar{c}(L_2), \label{Dominance_c} \\
d(L_1) &= d(L_2), \label{Dominance_d} \\
\mathcal{S}(L_1) &\subseteq \mathcal{S}(L_2), \label{Dominance_e} \\
(l_m^k(L_1))_{|\mathcal{M}_k|} &\le (l_m^k(L_2))_{|\mathcal{M}_k|}, \label{Dominance_f} \\
(O_m^k(L_1))_{|\mathcal{M}_k|} &\subseteq (O_m^k(L_2))_{|\mathcal{M}_k|}, \label{Dominance_g}\\
g(L_1) &\le g(L_2), \label{Dominance_h}\\
o(L_1) &\le o(L_2) \label{Dominance_i}
\end{align}
\end{subequations}
} 

\begin{prop}
Conditions (\ref{Dominance}) form a valid dominance rule.
\end{prop}

\begin{proof}
Labels are only compared if they belong to the same day~\eqref{Dominance_d}.
This implies that any further extension of the label accumulates resources in the same way (in particular, the constant $\boldsymbol{B}$ is identical in \eqref{REF_Time}) and is subject to the same bounds (daily working time limit $\boldsymbol{T}$ in \eqref{daily time limit} and daily driving limit $\boldsymbol{Dist}$ in \eqref{daily drive limit}).
As the REFs (\ref{REF_Time})-(\ref{REF_Drive_Day}) are all non-decreasing, the dominance rules are valid \citep{desaulniers1998unified}.
\end{proof}

\noindent \textbf{Inter-day Dominance.} The same-day dominance rule can be extended to compare labels across different days, allowing more aggressive pruning. Suppose two labels $L_1$ and $L_2$ satisfy the first two conditions ($i(L_1) = i(L_2)$ and $t(L_1) \le t(L_2)$), but $d(L_1) < d(L_2)$. Any feasible extension from $L_2$ can also be performed from $L_1$ if $L_1$ has sufficient time to catch up on a break. This break resets the daily resources $g(L)$ and $o(L)$, so no comparison is necessary for them. The only effect is an increase in the accumulated time by the break duration $\boldsymbol{B}$. The inter-day dominance conditions across days are therefore:

{\allowdisplaybreaks
\begin{subequations}\label{inter-day dominance}
\begin{align}
i(L_1) &= i(L_2), \label{inter-day dominance_a} \\
t(L_1) &\le t(L_2)-\boldsymbol{B}, \label{inter-day dominance_b} \\
\Bar{c}(L_1) &\le \Bar{c}(L_2), \label{inter-day dominance_c} \\
d(L_1) &< d(L_2), \label{inter-day dominance_d} \\
\mathcal{S}(L_1) &\subseteq \mathcal{S}(L_2), \label{inter-day dominance_e} \\
(l_m^k(L_1))_{|\mathcal{M}_k|} &\le (l_m^k(L_2))_{|\mathcal{M}_k|}, \label{inter-day dominance_f} \\
(O_m^k(L_1))_{|\mathcal{M}_k|} &\subseteq (O_m^k(L_2))_{|\mathcal{M}_k|} \label{inter-day dominance_g}
\end{align}
\end{subequations}
}

\noindent Note that the stronger condition $t(L_1) \le t(L_2)-\boldsymbol{B}$ is necessary, and $d(L_1) < d(L_2)$ alone is not sufficient.
E.g., consider a label $L_1$ late on day 1 and a label $L_2$ early on day 2.
An extension to another early customer on day 2 may be feasible for $L_2$ but infeasible for $L_1$, which may arrive too late after taking a necessary break.
Condition~\eqref{inter-day dominance_a} ensures there is sufficient time to take that break.

\begin{prop}
Conditions (\ref{inter-day dominance}) form a valid dominance rule.
\end{prop}
\begin{proof}
Let $L_1$ and $L_2$ be two labels satisfying (\ref{inter-day dominance}). Consider any feasible extension of $L_2$ along an arc $(i,j)$. Since $i(L_1)=i(L_2)$ and $t(L_1) + \boldsymbol{B} \le t(L_2)$, we can take an immediate break of duration $\boldsymbol{B}$ at $L_1$ before extending. Extensions are allowed to skip days if necessary, so this break is sufficient to catch up to $L_2$.
The break also resets daily-limited resources $g$ and $o$ to zero, ensuring that any extension feasible from $L_2$ is also feasible from $L_1$. Moreover, all other resources are no worse for $L_1$, and the time window constraints remain satisfied. Therefore, discarding $L_2$ does not eliminate any optimal solution, proving that the rule is valid.
\end{proof}

This enhancement allows $L_1$ to dominate $L_2$ even if they belong to different days with no requirement on $g(L_1)$ or $g(L_2)$ and $o(L_1)$ or $o(L_2)$, provided that an additional break can be taken to match the daily resource state of $L_2$. It preserves correctness while enabling more effective pruning of suboptimal labels, improving the efficiency of the labeling algorithm in the presence of driver breaks and daily limits.

The importance of Conditions~\eqref{inter-day dominance} can also be seen from the following.
For each label, two extensions are considered depending on whether a break is taken.
This has the potential to exponentially increase the number of non-dominated labels over each path.
With Conditions~\eqref{inter-day dominance}, however, unnecessary breaks that result in later days can typically be detected and eliminated, only maintaining both extensions when the decision is not obvious (e.g., taking a break now instead of making an additional delivery would facilitate an early delivery in the morning).

\subsection{Rolling-Space Branch-and-Price Algorithm}\label{RollingSpaceHorizonAlgorithm}

B\&P algorithms are highly effective for solving VRPs, but they often become computationally intractable for large-scale instances with complex constraints. To address this, this work proposes a \emph{Rolling-Space Branch-and-Price} (RS-B\&P) algorithm (see Algorithm \ref{alg:RSBP}) that integrates two key innovations: (i) a rolling-space clustering strategy that partitions the customers into overlapping spatial clusters to enable scalable and parallelizable route generation, and (ii) a multi-compartment feasibility check that efficiently enforces compatibility constraints without excessively increasing the size of the pricing problem.

\begin{figure}[t!]
    \centering
    \includegraphics[width=0.9\linewidth]{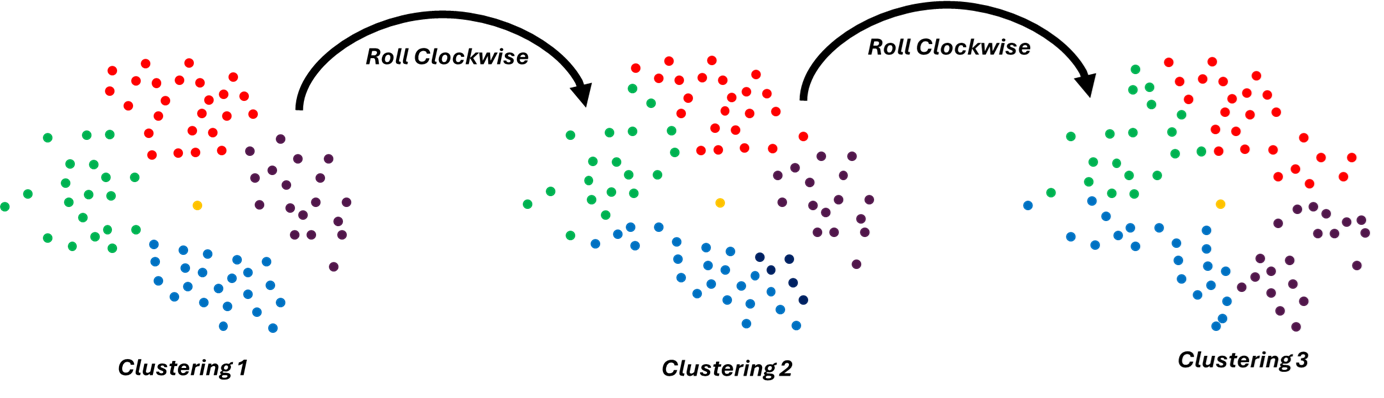}
    \caption{Rolling Space Example}
    \label{fig:Rolling Space Example}
\end{figure}

In the heuristic RS-B\&P, the set of customers is first partitioned into small clusters. The generation of overlapping clusters uses a sweep of a fixed-size window over the service area clockwise starting from the depot, as illustrated in Figure~\ref{fig:Rolling Space Example}. This ensures that each customer is considered in multiple clusters of close customers, improving the quality of candidate routes while keeping cluster-level pricing problems tractable. The exact B\&P algorithm described in Section~\ref{Exact B&P} is applied within each cluster, generating routes that are subsequently verified through a multi-compartment feasibility check. In other words, the generated routes initially ignore multi-compartment constraints (i.e., same-day dominance without conditions \eqref{Dominance_f}–\eqref{Dominance_g}, and inter-day dominance without conditions \eqref{inter-day dominance_f}–\eqref{inter-day dominance_g}). The multi-compartment feasibility check is only applied to routes that are likely to violate these constraints, specifically those approaching vehicle capacity or containing incompatible items. This corresponds to enforcing conditions \eqref{Dominance_f}–\eqref{Dominance_g} for same-day dominance and conditions \eqref{inter-day dominance_f}–\eqref{inter-day dominance_g} for inter-day dominance. This compatibility check is coupled with the rolling-space clustering strategy, ensuring that only promising routes from each cluster are retained for integration into the global RMP.

\begin{algorithm}[t!]
\footnotesize
\setstretch{1}
\caption{RS-B\&P Algorithm}
\label{alg:RSBP}
\textbf{Step 1: Initialization}\\
Construct initial rolling-space clusters\\
Generate initial feasible routes\\
Initialize global column pool with initial columns\\
\BlankLine
\textbf{Step 2: Root Node Processing}\\
Create root node $N_0$\\
Initialize the root RMP$_{N_0}$ with the columns' pool\\
\Repeat{no cluster produces a negative reduced-cost column}{
    \For{each cluster \textbf{in parallel}}{
        Solve pricing problem using dual prices from RMP$_{N_0}$\\
        \For{each generated candidate route $r$}{
            Perform feasibility check for multi-compartment constraints \\
            \If{$r$ is feasible and has negative reduced cost}{
                Add $r$ to RMP$_{N_0}$ and to the global pool \\
            }
        }
    }
    Re-solve RMP$_{N_0}$ to update dual prices\\
}
\If{RMP$_{N_0}$ is integer}{
    Store incumbent solution and terminate\\
}
\Else{
    Branch on a fractional arc to create child nodes
}
\BlankLine
\textbf{Step 3: Branch-and-Price}\\
\While{unexplored nodes exist}{
    Select next node $N$ according to node selection strategy (e.g., best-bound)\;
    Add arc branching constraints to the child RMP$_N$\\
    Initialize RMP$_N$ with all columns in the pool satisfying the branching constraints\\
    \Repeat{no cluster produces a negative reduced-cost column}{
        \For{each cluster \textbf{in parallel}}{
            Solve pricing with duals from RMP$_N$\\
            \For{each generated candidate route $r$}{
            Perform feasibility check for multi-compartment constraints \\
            \If{$r$ is feasible and has negative reduced cost}{
                Add $r$ to RMP$_{N}$ and to global pool \\
            }
            }
        }
        Re-solve RMP$_N$\;
    }
    \If{RMP$_N$ is integer}{
        Update incumbent if better; fathom node
    }
    \Else{
        Branch on a fractional arc to create child nodes
    }
}
\BlankLine
\textbf{Step 4: Return incumbent solution.}
\end{algorithm}

The verification of the multi-compartment constraints for all candidate routes whose total load is close to the vehicle capacity (within a predefined threshold $\kappa$) is achieved by solving the following feasibility problem for each such route $r \in \Omega$ (associated with vehicle $k \in \mathcal{K}$).

{\small
\allowdisplaybreaks
\centering
\begin{align}
\label{Feasibility Model}
\min_{w,r} \hspace{2mm} & 0, \\[2mm]
\label{eq:demand} 
\text{s.t.} \hspace{2mm} 
& \sum_{m \in \mathcal{M}_k} r_{nm} = 1, 
&& \forall n \in \mathcal{N}_r, \\
\label{eq:compartmentcompatability} 
& r_{nm} \leq \boldsymbol{b}_{nm}, 
&& \forall n \in \mathcal{N}_r,\, m \in \mathcal{M}_k, \\
\label{eq:capacity} 
& \sum_{n \in \mathcal{N}_r} \boldsymbol{Q}_n r_{nm} \leq \boldsymbol{L}_{mk}, 
&& \forall m \in \mathcal{M}_k, \\
\label{eq:incompatibility} 
& r_{im} + r_{jm} \leq 1 + \boldsymbol{f}_{ij}, 
&& \forall i,j \in \mathcal{N}_r,\, m \in \mathcal{M}_k, \\
\label{eq:fractional} 
& r_{nm} \in \{0,1\} 
&& \forall n \in \mathcal{N}_r,\, m \in \mathcal{M}_k.
\end{align}
}

\noindent where $r_{nm} \in \{0,1\}$ is 1 if customer $n$ is assigned to compartment $m$, 0 otherwise. Constraints (\ref{eq:demand}) ensure that a customer is served by exactly one compartment. Constraints (\ref{eq:compartmentcompatability}) ensure compatibility between the customers and the compartments. Constraints (\ref{eq:capacity}) impose a capacity limit on each compartment to ensure that the total load assigned does not exceed the compartment's maximum capacity. Constraints (\ref{eq:incompatibility}) ensure incompatibility restrictions by preventing any pair of incompatible customers from being assigned to the same compartment. Constraints (\ref{eq:fractional}) ensure that the restrictions on the variables are met. The routes that are not feasible from a loading perspective are not added to the RMP.

By solving each pricing problem at the cluster level, RS-B\&P exploits both spatial decomposition and parallelism to accelerate route generation. The resulting routes are aggregated into a single master problem, while the use of overlapping clusters enhances solution quality by ensuring broader route coverage. Multi-compartment constraints are enforced efficiently through a dedicated post-processing step. Only candidate routes whose total load is within a predefined threshold $\kappa$ of the vehicle capacity, or that include customers assigned to incompatible product compartments, are validated through the feasibility model above. Routes failing this test are discarded before being added to the RMP. This selective validation strategy prevents unnecessary growth of the pricing problem, allowing the algorithm to scale effectively to large MCVRPMTW instances.

\subsection{Acceleration Strategies}

This section proposes several acceleration techniques for the B\&P algorithm. The first two methods that enable the acceleration of the pricing problem. The third method accelerates the master problem. The fourth method enhances the branching process. These techniques have been adopted from the literature.

\subsubsection{Pricing Problem Acceleration}

Arc filtering is employed to reduce the number of arcs considered during label extension. This technique preemptively eliminates arcs that cannot lead to feasible or cost-effective routes based on vehicle capacity, time windows, or dominance conditions discussed above. It is also employed to break the symmetry that can be obtained from the multiple time windows or multi-compartment attributes. If one of the conditions is not met, the arc is deemed infeasible; otherwise, the extension is permitted. To break symmetry, the algorithm selects the arc(s) with the earliest possible time window and the smallest index possible compartment.

The algorithms also implement an aggressive dominance rule. The aggressive same-day dominance \citep{costa2019exact} defines a relaxed dominance criterion between two labels by incorporating $\epsilon$ tolerances on resource values. The procedure first verifies that both labels terminate at the same node and that the first label has visited a subset of the customers visited by the second. It then permits the first label to exceed the second in cost, time, or load by small predefined tolerance values ($\epsilon$), thereby allowing limited degradations across resources in exchange for broader label pruning. Early iterations of the column generation process apply the aggressive dominance rule to reduce the number of propagated labels and accelerate convergence. Once the reduced costs stabilize, the algorithm switches to the exact Pareto dominance rule. The same logic applies to the aggressive inter-day dominance. The aggressive same-day dominance is as follows

{\small
\allowdisplaybreaks
\begin{equation}\label{EpsilonDominance}
    \begin{cases}
      i(L_1)=i(L_2), & \\
      t(L_1) \leq t(L_2) + \boldsymbol{\varepsilon_t}, & \\
      \Bar{c}(L_1) \leq \Bar{c}(L_2) + \boldsymbol{\varepsilon_c}, & \\
      d(L_1) = d(L_2), & \\
      \mathcal{S}(L_1) \subseteq \mathcal{S}(L_2), & \\
      (l_m^k(L_1))_{|\mathcal{M}_k|} \leq (l_m^k(L_2))_{|\mathcal{M}_k|}, & \\
      (O_m^k(L_1))_{|\mathcal{M}_k|} \subseteq (O_m^k(L_2))_{|\mathcal{M}_k|} & \\
      g(L_1) \leq g(L_2) + \boldsymbol{\varepsilon_g}, & \\
      o(L_1) \leq o(L_2) + \boldsymbol{\varepsilon_o}, & \\
    \end{cases}       
\end{equation}
}

\subsubsection{Master Problem Acceleration}

 Rather than adding all columns (routes) with negative reduced costs to the RMP, the column selection prioritizes a subset based on the reduced costs, i.e., inserting the columns with the most negative reduced costs. The selection function first sorts the input routes in ascending order of reduced cost, prioritizing those with the most negative values. Then, it iterates through the sorted list, adding up to a selected maximum number of routes (or percentage threshold) to the master problem. A second function selects high-quality, mutually disjoint routes for the master problem. It first sorts the routes by reduced cost, then iterates through them to identify those that are both cost-effective and do not overlap with already selected routes. Selected routes are added up to a specified maximum, and customer indices in each accepted route are marked to prevent overlap. 

\subsubsection{Branch \& Bound Acceleration}

Primal diving in the VRPMCs \citep{bani2024combining} is an acceleration technique used within the B\&B framework to rapidly identify high-quality feasible solutions, thereby enhancing pruning efficiency and guiding the search toward promising regions of the solution space. It operates by simulating a depth-first traversal of the B\&B tree, fixing variables based on fractional values from the current LP relaxation, typically favoring those closest to integer values or most impactful on the objective. By prioritizing rapid construction of feasible primal solutions, the heuristic can improve upper bounds early in the process, reduce the number of nodes explored, and accelerate overall convergence.

\section{Computational Results}\label{section:5}

This section reports extensive computational experiments conducted to evaluate the performance of the B\&P and RS-B\&P algorithms and to discuss managerial findings on key characteristics of the MCVRPMTW. The algorithm was implemented in C++ and compiled as 64-bit single-threaded code using MS Visual Studio 2019. GUROBI 12.0.2, with default parameters, was used to solve the RMPs. Experiments were run on a server equipped with dual-socket Intel Xeon Gold 6226 CPUs, each featuring 24 cores at 2.7~GHz. The performance of a single thread on this system is comparable to that of a typical desktop processor. A maximum runtime of 36,000 seconds was applied to all computational runs.

The remainder of this section is structured as follows. Section~\ref{section5:Instances} describes the benchmark instances and specifies the parameter settings for the studied attributes. Section~\ref{section5:ResultsSummary} provides a summary of the algorithm’s performance. Sections~\ref{section5:AccelerationStrategies} and \ref{section5:ImpactAttributes} analyze the effects of acceleration strategies and the impact of different MCVRPMTW attributes, respectively. Finally, Section~\ref{section5:ManagerialInsights} presents the managerial insights derived from the results.

\subsection{Instances and Attribute Settings}\label{section5:Instances}

The instances are inspired by the real-world application that motivated this research. In this context, the company receives weekly orders from various US States to fulfill. Vehicles load weekly orders, deliver them over the course of a week, and return to the depot for the next cycle. Each day, each customer has a single time window during which the delivery is possible. The goal is to deliver the weekly orders while minimizing the total distance and the number of vehicles used. Drivers must take a 10-hour break at the end of each daily shift and cannot exceed 580 miles or work more than 12 hours per day. All these practical constraints are enforced during route generation and feasibility checks, ensuring that the computed solutions are directly implementable in the company's operational setting. 

From this real-world application, a set of 14 realistic instances were generated and classified into two classes: (1) Instances with up to 40 customers solved using our B\&P algorithm, (2) Instances with up to 400 customers solved using B\&P enhanced with clustering as described in Section \ref{RollingSpaceHorizonAlgorithm}. A high artificial cost is added to each route $r \in \Omega$ to ensure the minimization of the number of vehicles. The instances consider vehicles with six compartments (all of equal capacity), with a total capacity of 40,000 for both instance classes. The instances also include eight item categories (or compatibility numbers) to model item-to-item and item-to-compartment compatibility. Each item is assigned to a category based on the company data for each load. For instance, certain orders must be delivered from the rear and are therefore restricted to rear compartments, while some loads cannot be placed on top of others. The 6-compartment case is detailed in Appendix~\ref{C}. Unless otherwise specified, all computational results refer to the 6-compartment configuration.

\subsection{Summary of Computational Results}\label{section5:ResultsSummary}

The B\&P and RS-B\&P algorithms were applied on the first classes of instances (M1 to M7). Table \ref{tab:BP Results on M Instances} shows computational results and reports the number of customers (\#Cust),  total load (TL), objective value (Obj), which is the best solution found, number of vehicles used (\#Veh), number of columns (\#Col), number of nodes (\#Nod), and total computation time (Time). For B\&P, we report whether the instance was solved to optimality (Opt), and for RS-B\&P (with acceleration strategies), we report the gap relative to B\&P (Gap).

\begin{table}[t!]
  \centering
  \caption{Results on M Instances}
  \label{tab:BP Results on M Instances}%
  \scalebox{0.65}{
    \begin{tabular}{c|c|c|cccccc|cccccc}
    \toprule
    \multirow{2}[4]{*}{\textbf{\#Inst}} & \multirow{2}[4]{*}{\textbf{\#Cust}} & \multirow{2}[4]{*}{\textbf{TL}} & \multicolumn{6}{c|}{\textbf{B\&P}}       & \multicolumn{6}{c}{\textbf{RS-B\&P}} \\
\cmidrule{4-15}          &       &       & \textbf{Obj} & \textbf{\#Veh} & \textbf{\#Col} & \textbf{\#Nod} & \textbf{Time} & \textbf{Opt} & \textbf{Obj} & \textbf{\#Veh} & \textbf{\#Col} & \textbf{\#Nod} & \textbf{Time} & \textbf{Gap} \\
    \midrule
    \textbf{M1} & 10    & 21{,}292.82 & 1{,}362.80 & 1     & 18{,}863 & 1     & 367.13 & \ding{51}    & 1{,}363.71 & 1     & 5{,}491  & 1     & 3.47  & 0.07\% \\
    \textbf{M2} & 15    & 33{,}835.52 & 1{,}582.45 & 1     & 29{,}536 & 6     & 2{,}520.99 & \ding{51}    & 1{,}592.84 & 1     & 8{,}456  & 3     & 39.21 & 0.65\% \\
    \textbf{M3} & 20    & 47{,}658.39 & 2{,}718.16 & 2     & 5{,}8716 & 42    & 5{,}172.29 & \ding{51}    & 2{,}795.40 & 2     & 13{,}022 & 19    & 64.26 & 2.76\% \\
    \textbf{M4} & 25    & 59{,}608.06 & 2{,}935.61 & 2     & 64{,}200 & 80    & 8{,}060.41 & \ding{51}    & 3{,}024.30 & 2     & 2{,}0054 & 37    & 91.54 & 2.93\% \\
    \textbf{M5} & 30    & 71{,}263.73 & 3{,}153.05 & 2     & 77{,}580 & 158   & 16{,}120.82 & \ding{51}    & 3{,}275.16 & 2     & 30{,}884 & 73    & 167.68 & 3.73\% \\
    \textbf{M6} & 35    & 84{,}511.10 & 4{,}294.69 & 3     & 88{,}252 & 289   & 36{,}000.00 & \ding{55}    & 4{,}530.97 & 3     & 47{,}561 & 133   & 389.47 & 5.21\% \\
    \textbf{M7} & 40    & 94{,}936.14 & 4{,}842.33 & 3     & 11{,}7432 & 806   & 36{,}000.00 & \ding{55}    & 4{,}967.34 & 3     & 73{,}243 & 370   & 557.16 & 2.52\% \\
    \midrule
    \textbf{Avg} & \textbf{25} & \textbf{59{,}015.11} & \textbf{2{,}984.16} & \textbf{2} & \textbf{64{,}940} & \textbf{198} & \textbf{14{,}891.66} & \textbf{-} & \textbf{3{,}078.53} & \textbf{2} & \textbf{25{,}807} & \textbf{91} & \textbf{187.54} & \textbf{2.55\%} \\
    \bottomrule
    \end{tabular}%
}
\end{table}%

The results highlight that, as the number of customers increases, the computation times grow significantly. The number of vehicles used also scales proportionally with the customer count. Overall, the B\&P algorithm reaches optimal solutions on the first five instances, with an average solving time of approximately 14{,}892 seconds. On the other hand, the RS-B\&P algorithm achieves comparable solutions, with an average gap of 2.55\% in an average time of 187.54 seconds. It highlights that the RS-B\&P algorithm produces high-quality solutions close to the B\&P algorithm, while significantly reducing the runtime.

For the large instances (L1 to L7), Table \ref{tab:Comparison between BP and LH for L Instances} compares the performance of the proposed RS-B\&P algorithm to a labeling-based heuristic (LH). The LH constructs routes using a labeling procedure. Each partial route (label) is iteratively extended by evaluating all feasible next customers that respect time windows, multi-compartment constraints, and driver regulations. A route grows until the vehicle reaches its maximum load of 40,000, or no additional customers can be feasibly added. Once a route is completed, it is fixed, and the procedure proceeds to generate new routes for the remaining unassigned customers. This process continues until all customers are assigned to routes. Conceptually, the heuristic resembles a simplified column-generation approach, as it generates feasible routes on the fly, prioritizing efficient vehicle utilization while maintaining feasibility.

\begin{table}[t!]
\centering
\scalebox{0.7}{
\caption{Comparison between RS-B\&P and LH for L Instances}
\label{tab:Comparison between BP and LH for L Instances}%
    \begin{tabular}{ccc|cccc|cccc}
    \toprule
    \multirow{2}[4]{*}{\textbf{\#Inst}} & \multirow{2}[4]{*}{\textbf{\#Cust}} & \multirow{2}[4]{*}{\textbf{TL}} & \multicolumn{4}{c|}{\textbf{RS-B\&P}} & \multicolumn{4}{c}{\textbf{LH}} \\
\cmidrule{4-11}          &       &       & \textbf{Obj} & \textbf{\#Veh} & \textbf{\#Col} & \textbf{Time} & \textbf{Obj} & \textbf{\#Veh} & \textbf{\#Col} & \textbf{Time} \\
    \midrule
    \textbf{L1} & 50    & 117{,}049.69 & 5{,}346.12 & 3     & 16{,}473 & 1{,}259.27 & 6{,}952.33 & 4     & 7{,}413  & 0.55 \\
    \textbf{L2} & 60    & 139{,}881.54 & 7{,}082.29 & 4     & 25{,}368 & 1{,}385.20 & 8{,}528.56 & 5     & 11{,}416 & 0.55 \\
    \textbf{L3} & 75    & 200{,}675.54 & 12{,}401.15 & 6     & 39{,}067 & 1{,}662.24 & 12{,}646.30 & 8     & 17{,}580 & 0.45 \\
    \textbf{L4} & 100   & 259{,}114.57 & 14{,}693.57 & 7     & 60{,}163 & 2{,}327.13 & 15{,}827.30 & 9     & 27{,}073 & 0.71 \\
    \textbf{L5} & 200   & 467{,}837.82 & 28{,}438.15 & 13    & 92{,}651 & 3{,}723.41 & 29{,}440.10 & 16    & 41{,}693 & 1.01 \\
    \textbf{L6} & 300   & 725{,}175.74 & 41{,}733.68 & 20    & 142{,}682 & 4{,}468.09 & 43{,}856.80 & 23    & 64{,}207 & 1.08 \\
    \textbf{L7} & 400   & 953{,}081.54 & 55{,}556.30 & 25    & 219{,}730 & 6{,}947.47 & 57{,}015.80 & 28    & 98{,}879 & 1.06 \\
    \midrule
    \textbf{Avg} & \textbf{169} & \textbf{408{,}973.78} & \textbf{23{,}607.32} & \textbf{11} & \textbf{85{,}162} & \textbf{3{,}110.40} & \textbf{24{,}895.31} & \textbf{13} & \textbf{38{,}323} & \textbf{0.77} \\
    \bottomrule
    \end{tabular}%
}
\end{table}%

As the problem size increases, the solution times grow reasonably, with computation time reaching over 6{,}900 seconds for the largest instance (L7). The number of vehicles scales with customer counts, from 3 to 25. The RS-B\&P algorithm remains effective, maintaining reasonable runtimes and demonstrating scalability up to 400 customers. The LH heuristic is very fast, solving instances in under 1 second on average, and produces feasible solutions of reasonable quality. However, as shown in Table~\ref{tab:Comparison between BP and LH for L Instances}, the RS-B\&P algorithm consistently achieves higher-quality solutions, reducing total cost by roughly 5\% on average. More importantly, RS-B\&P consistently uses fewer vehicles, which is a critical operational metric in practice.

As shown in Table~\ref{tab:BP Results on M Instances}, the RS-B\&P is fast and delivers high-quality solutions on small benchmark instances. Furthermore, on large benchmark instances, Table~\ref{tab:Comparison between BP and LH for L Instances} highlights the trade-off between computational efficiency and solution quality. The LH heuristic provides very fast \emph{good-enough} solutions. At the same time, the RS-B\&P achieves superior performance, particularly in more constrained or large-scale settings. The RS-B\&P also scales efficiently to real-life, large-scale instances, where it identifies substantial operational savings compared to the LH. These findings, discussed in detail in Section~\ref{section5:ManagerialInsights}, illustrate both the practical relevance and the robustness of the proposed method in complex logistics settings.

\subsection{Impact of Acceleration Strategies}\label{section5:AccelerationStrategies}

Table \ref{tab:Impact of Acceleration Strategies} presents the impact of disabling acceleration strategies on the performance of the RS-B\&P algorithm across both medium- and large-scale instances. These strategies include Warmstart (W), Arc Filtering (AF), Aggressive/Pareto Dominance (APD), Column Selection (CS), and Primal Diving Heuristic (PDH). For each instance, performance is compared under five configurations where one strategy is removed at a time (No W, No AF \& APD, No CS, No PDH), and benchmarked against the RS-B\&P version where all strategies are enabled. Each configuration corresponds to a step in the solving process, i.e., pricing problem, master problem, and branching. We report the objective value, number of vehicles, and computation time.

\begin{table}[t!]
\centering
\scalebox{0.5}{
  \caption{Impact of Acceleration Strategies}
  \label{tab:Impact of Acceleration Strategies}%
    \begin{tabular}{c|c|c|ccc|ccc|ccc|ccc|ccc}
    \toprule
    \multirow{2}[4]{*}{\textbf{\#Inst}} & \multirow{2}[4]{*}{\textbf{\#Cust}} & \multirow{2}[4]{*}{\textbf{TL}} & \multicolumn{3}{c|}{\textbf{No W}} & \multicolumn{3}{c|}{\textbf{No AF \& A/PD}} & \multicolumn{3}{c|}{\textbf{No CS}} & \multicolumn{3}{c|}{\textbf{No PDH}} & \multicolumn{3}{c}{\textbf{RS-B\&P}} \\
\cmidrule{4-18}          &       &       & \textbf{Obj} & \textbf{\#Veh} & \textbf{Time} & \textbf{Obj} & \textbf{\#Veh} & \textbf{Time} & \textbf{Obj} & \textbf{\#Veh} & \textbf{Time} & \textbf{Obj} & \textbf{\#Veh} & \textbf{Time} & \textbf{Obj} & \textbf{\#Veh} & \textbf{Time} \\
    \midrule
    \textbf{M1} & 10    & 21{,}292.82 & 1{,}363.71 & 1     & 3.92  & 1{,}363.71 & 1     & 7.63  & 1{,}363.71 & 1     & 4.86  & 1{,}363.71 & 1     & 9.71  & 1{,}363.71 & 1     & 3.47 \\
    \textbf{M2} & 15    & 33{,}835.52 & 1{,}592.84 & 1     & 44.31 & 1{,}592.84 & 1     & 86.27 & 1{,}592.84 & 1     & 54.90 & 1{,}592.84 & 1     & 109.80 & 1{,}592.84 & 1     & 39.21 \\
    \textbf{M3} & 20    & 47{,}658.39 & 2{,}795.40 & 2     & 72.61 & 2{,}795.40 & 2     & 141.37 & 2{,}795.40 & 2     & 89.96 & 2{,}795.40 & 2     & 179.93 & 2{,}795.40 & 2     & 64.26 \\
    \textbf{M4} & 25    & 59{,}608.06 & 3{,}024.30 & 2     & 103.44 & 3{,}024.30 & 2     & 201.39 & 3{,}024.30 & 2     & 128.16 & 3{,}024.30 & 2     & 256.32 & 3{,}024.30 & 2     & 91.54 \\
    \textbf{M5} & 30    & 71{,}263.73 & 3{,}275.16 & 2     & 189.48 & 3{,}275.16 & 2     & 368.90 & 3{,}275.16 & 2     & 234.76 & 3{,}275.16 & 2     & 469.51 & 3{,}275.16 & 2     & 167.68 \\
    \textbf{M6} & 35    & 84{,}511.10 & 4{,}530.97 & 3     & 440.11 & 4{,}530.97 & 3     & 856.84 & 4{,}530.97 & 3     & 545.26 & 4{,}530.97 & 3     & 1{,}090.53 & 4{,}530.97 & 3     & 389.47 \\
    \textbf{M7} & 40    & 94{,}936.14 & 4{,}967.34 & 3     & 629.59 & 4{,}967.34 & 3     & 1{,}225.75 & 4{,}967.34 & 3     & 780.02 & 4{,}967.34 & 3     & 1{,}560.04 & 4{,}967.34 & 3     & 557.16 \\
    \midrule
    \textbf{Avg} & \textbf{25} & \textbf{59{,}015.11} & \textbf{3{,}078.53} & \textbf{2} & \textbf{211.92} & \textbf{3078.53} & \textbf{2} & \textbf{412.59} & \textbf{3{,}078.53} & \textbf{2} & \textbf{262.56} & \textbf{3{,}078.53} & \textbf{2} & \textbf{525.12} & \textbf{3{,}078.53} & \textbf{2} & \textbf{187.54} \\
    \midrule
    \textbf{L1} & 50    & 117{,}049.69 & 5{,}346.12 & 3     & 1{,}422.97 & 5{,}385.32 & 3     & 2{,}770.39 & 5{,}385.32 & 3     & 1{,}762.98 & 5{,}385.32 & 3     & 3{,}540.53 & 5{,}346.12 & 3     & 1{,}259.27 \\
    \textbf{L2} & 60    & 139{,}881.54 & 7{,}082.29 & 4     & 1{,}565.27 & 7{,}126.75 & 4     & 3{,}047.43 & 7{,}126.75 & 4     & 1{,}939.27 & 7{,}126.75 & 4     & 4{,}851.63 & 7{,}082.29 & 4     & 1{,}385.20 \\
    \textbf{L3} & 75    & 200{,}675.54 & 12{,}401.15 & 6     & 1{,}878.33 & 12{,}481.60 & 6     & 3{,}656.92 & 12{,}481.60 & 6     & 2{,}327.13 & 17{,}361.61 & 9     & 7{,}200.00 & 12{,}401.15 & 6     & 1{,}662.24 \\
    \textbf{L4} & 100   & 259{,}114.57 & 14{,}693.57 & 7     & 2{,}629.66 & 14{,}910.70 & 7     & 5{,}119.68 & 14{,}910.70 & 7     & 3{,}257.98 & 20{,}571.00 & 10    & 7{,}200.00 & 14{,}693.57 & 7     & 2{,}327.13 \\
    \textbf{L5} & 200   & 467{,}837.82 & 28{,}438.15 & 13    & 4{,}207.45 & 37{,}422.84 & 15    & 7{,}200.00 & 28{,}786.80 & 13    & 5{,}212.77 & 39{,}813.41 & 19    & 7{,}200.00 & 28{,}438.15 & 13    & 3{,}723.41 \\
    \textbf{L6} & 300   & 72{,}5175.74 & 41{,}733.68 & 20    & 5{,}048.94 & 50{,}485.80 & 21    & 7{,}200.00 & 42{,}071.50 & 20    & 6{,}255.32 & 58{,}427.15 & 28    & 7{,}200.00 & 41{,}733.68 & 20    & 4{,}468.09 \\
    \textbf{L7} & 400   & 953{,}081.54 & 55{,}556.30 & 25    & 7{,}850.64 & 67{,}218.96 & 28    & 7{,}200.00 & 58{,}015.80 & 27    & 7{,}200.00 & 77{,}778.82 & 35    & 7{,}200.00 & 55{,}556.30 & 25    & 6{,}947.47 \\
    \midrule
    \textbf{Avg} & \textbf{169} & \textbf{408{,}973.78} & \textbf{23{,}607.32} & \textbf{11} & \textbf{3{,}514.75} & \textbf{27{,}861.71} & \textbf{12} & \textbf{5{,}170.63} & \textbf{24{,}111.21} & \textbf{11} & \textbf{3{,}993.64} & \textbf{32{,}352.01} & \textbf{15} & \textbf{6{,}341.74} & \textbf{23{,}607.32} & \textbf{11} & \textbf{3{,}110.40} \\
    \bottomrule
    \end{tabular}%
}
\end{table}%

The results demonstrate that all five strategies contribute significantly to improving computational efficiency, particularly for large instances. Removing warmstart (No W) or column selection (No CS) leads to moderate increases in runtime, while removing arc filtering and aggressive dominance (No AF \& AD) has a stronger effect, more than doubling the average execution time. However, the absence of the primal diving heuristic (No PDH) is especially detrimental. In particular, for large instances, it results in the highest computation times and dramatically increases the objective (up to 27\%). In some cases, the solver reaches the time limit without improving the solution. Overall, the RS-B\&P configuration outperforms all alternatives, highlighting that combining all strategies yields faster convergence and higher solution quality.

\subsection{Impact of the MCVRPMTW Attributes} \label{section5:ImpactAttributes}

An extensive sensitivity analysis was performed to evaluate various parameters of our RS-B\&P algorithm, including the impact of compartment attributes, the impact of time windows, and the impact of the clustering approach.

\subsubsection{Impact of Compartment Attributes}

\begin{table}[t!]
\centering
\scalebox{0.6}{
  \caption{Impact of Compartment Attributes}
  \label{tab:Impact of Compartment Attributes}%
    \begin{tabular}{c|c|c|ccc|ccc|ccc|ccc}
    \toprule
    \multirow{2}[4]{*}{\textbf{\#Inst}} & \multirow{2}[4]{*}{\textbf{\#Cust}} & \multirow{2}[4]{*}{\textbf{TL}} & \multicolumn{3}{c|}{\textbf{2 Compartments}} & \multicolumn{3}{c|}{\textbf{4 Compartments}} & \multicolumn{3}{c|}{\textbf{6 Compartments}} & \multicolumn{3}{c}{\textbf{8 Compartments}} \\
\cmidrule{4-15}          &       &       & \textbf{Obj} & \textbf{\#Veh} & \textbf{Time} & \textbf{Obj} & \textbf{\#Veh} & \textbf{Time} & \textbf{Obj} & \textbf{\#Veh} & \textbf{Time} & \textbf{Obj} & \textbf{\#Veh} & \textbf{Time} \\
    \midrule
    \textbf{M1} & 10    & 21{,}292.82 & 1{,}500.08 & 2     & 2.67  & 1{,}431.90 & 1     & 2.89  & 1{,}363.71 & 1     & 3.47  & 1{,}363.71 & 1     & 3.81 \\
    \textbf{M2} & 15    & 33{,}835.52 & 1{,}752.12 & 6     & 30.16 & 1{,}672.48 & 5     & 32.68 & 1{,}592.84 & 1     & 39.21 & 1{,}592.84 & 1     & 43.14 \\
    \textbf{M3} & 20    & 47{,}658.39 & 3{,}074.94 & 3     & 49.43 & 2{,}935.17 & 3     & 53.55 & 2{,}795.40 & 2     & 64.26 & 2{,}795.40 & 2     & 70.69 \\
    \textbf{M4} & 25    & 59{,}608.06 & 3{,}326.73 & 3     & 70.42 & 3{,}175.52 & 3     & 76.28 & 3{,}024.30 & 2     & 91.54 & 3{,}024.30 & 2     & 100.70 \\
    \textbf{M5} & 30    & 71{,}263.73 & 3{,}602.68 & 3     & 128.99 & 3{,}438.92 & 3     & 139.74 & 3{,}275.16 & 2     & 167.68 & 3{,}275.16 & 2     & 184.45 \\
    \textbf{M6} & 35    & 84{,}511.10 & 4{,}984.06 & 5     & 299.60 & 4{,}757.51 & 4     & 324.56 & 4{,}530.97 & 3     & 389.47 & 4{,}530.97 & 3     & 428.42 \\
    \textbf{M7} & 40    & 94{,}936.14 & 5{,}464.07 & 5     & 428.58 & 5{,}215.70 & 4     & 464.30 & 4{,}967.34 & 3     & 557.16 & 4{,}967.34 & 3     & 612.87 \\
    \midrule
    \textbf{Avg} & \textbf{25} & \textbf{59{,}015.11} & \textbf{3{,}386.38} & \textbf{4} & \textbf{144.26} & \textbf{3{,}232.46} & \textbf{3} & \textbf{156.29} & \textbf{3{,}078.53} & \textbf{2} & \textbf{187.54} & \textbf{3{,}078.53} & \textbf{2} & \textbf{206.30} \\
    \midrule
    \textbf{L1} & 50    & 117{,}049.69 & 5{,}880.73 & 5     & 968.67 & 5{,}613.42 & 4     & 1{,}049.39 & 5{,}346.12 & 3     & 1{,}259.27 & 5{,}346.12 & 3     & 1{,}385.20 \\
    \textbf{L2} & 60    & 139{,}881.54 & 7{,}790.52 & 6     & 1{,}065.54 & 7{,}436.40 & 6     & 1{,}154.33 & 7{,}082.29 & 4     & 1{,}385.20 & 7{,}082.29 & 4     & 1{,}523.72 \\
    \textbf{L3} & 75    & 200{,}675.54 & 13{,}641.26 & 9     & 1{,}278.64 & 13{,}021.20 & 8     & 1{,}385.20 & 12{,}401.15 & 6     & 1{,}662.24 & 12{,}401.15 & 6     & 1{,}828.46 \\
    \textbf{L4} & 100   & 259{,}114.57 & 16{,}162.93 & 11    & 1{,}790.10 & 15{,}428.25 & 10    & 1{,}939.27 & 14{,}693.57 & 7     & 2{,}327.13 & 14{,}693.57 & 7     & 2{,}559.84 \\
    \textbf{L5} & 200   & 467{,}837.82 & 31{,}281.97 & 19    & 2{,}864.16 & 29{,}860.06 & 17    & 3{,}102.84 & 28{,}438.15 & 13    & 3{,}723.41 & 28{,}438.15 & 13    & 4{,}095.75 \\
    \textbf{L6} & 300   & 725{,}175.74 & 45{,}907.05 & 25    & 3{,}436.99 & 43{,}820.37 & 23    & 3{,}723.41 & 41{,}733.68 & 20    & 4{,}468.09 & 41{,}733.68 & 20    & 4{,}914.90 \\
    \textbf{L7} & 400   & 953{,}081.54 & 61{,}111.93 & 30    & 5{,}344.21 & 58{,}334.12 & 27    & 5{,}789.56 & 55{,}556.30 & 25    & 6{,}947.47 & 55{,}556.30 & 25    & 7{,}642.21 \\
    \midrule
    \textbf{Avg} & \textbf{169} & \textbf{408{,}973.78} & \textbf{25{,}968.05} & \textbf{15} & \textbf{2{,}392.61} & \textbf{24{,}787.69} & \textbf{14} & \textbf{2{,}592.00} & \textbf{23{,}607.32} & \textbf{11} & \textbf{3{,}110.40} & \textbf{23{,}607.32} & \textbf{11} & \textbf{3{,}421.44} \\
    \bottomrule
    \end{tabular}%
}
\end{table}%

Consider the impact of the number of compartments on the total distance, the number of vehicles, and the execution time. Table \ref{tab:Impact of Compartment Attributes} compares the performance of a vehicle routing problem with varying numbers of vehicle compartments (2, 4, 6, and 8) across the medium and large instances. Appendix \ref{C} summarizes the tested attribute settings. As the number of compartments increases, the overall trend shows a consistent reduction in both the objective value and the number of vehicles required. This improvement is due to the increased flexibility of loading diverse product types into separate compartments within a single vehicle. In the medium instance class, the average objective value decreases from 3386.38 to 3078.53, while the average number of vehicles drops from 4 to 2. A similar pattern holds for large instances, where the objective drops from 25{,}968.05 to 23{,}607.32, and vehicles decrease from 15 to 11. These results confirm that increasing the compartment count leads to more efficient routing solutions with better load consolidation.

\begin{figure}[t!]
    \centering
    \begin{subfigure}[b]{0.48\textwidth}
        \centering
        \begin{tikzpicture}
        \begin{axis}[
            width=\textwidth,
            xlabel={Number of Compartments},
            ylabel={Objective},
            ymin=10000, ymax=65000,
            yticklabel style={/pgf/number format/fixed},
            xtick={2,4,6,8},
            legend style={at={(0.5,-0.2)}, anchor=north, legend columns=2},
            ymajorgrids=true,
            grid style=dashed,
        ]
        \addplot[
            color=blue,
            mark=*,
            thick
        ] coordinates {
            (2,16162.93)
            (4,15428.25)
            (6,14693.57)
            (8,14693.57)
        };
        \addlegendentry{L4}
        \addplot[
            color=red,
            mark=*,
            thick
        ] coordinates {
            (2,31281.97)
            (4,29860.06)
            (6,28438.15)
            (8,28438.15)
        };
        \addlegendentry{L5}
        \addplot[
            color=brown,
            mark=*,
            thick
        ] coordinates {
            (2,45907.05)
            (4,43820.37)
            (6,41733.68)
            (8,41733.68)
        };
        \addlegendentry{L6}
        \addplot[
            color=orange,
            mark=*,
            thick
        ] coordinates {
            (2,61111.93)
            (4,58334.12)
            (6,55556.30)
            (8,55556.30)
        };
        \addlegendentry{L7}
        \end{axis}
        \end{tikzpicture}
        \caption{Objective Evolution with Compartments}
        \label{fig:L5_obj_gap}
    \end{subfigure}
    \hfill
    \begin{subfigure}[b]{0.48\textwidth}
        \centering
        \begin{tikzpicture}
        \begin{axis}[
            width=\textwidth,
            xlabel={Number of Compartments},
            ylabel={Execution Time},
            xtick={2,4,6,8},
            ymin=0, ymax=7700,
            legend style={at={(0.5,-0.2)}, anchor=north, legend columns=2},
            ymajorgrids=true,
            grid style=dashed,
        ]
        \addplot[
            color=blue,
            mark=diamond*,
            thick,
            dashed
        ] coordinates {
            (2, 1790.10)
            (4, 1939.27)
            (6, 2327.13)
            (8, 2559.84)
        };
        \addlegendentry{L4}
        \addplot[
            color=red,
            mark=diamond*,
            thick,
            dashed
        ] coordinates {
            (2, 2864.16)
            (4, 3102.84)
            (6, 3723.41)
            (8, 4095.75)
        };
        \addlegendentry{L5}
        \addplot[
            color=brown,
            mark=diamond*,
            thick,
            dashed
        ] coordinates {
            (2, 3436.99)
            (4, 3723.41)
            (6, 4468.09)
            (8, 4914.90)
        };
        \addlegendentry{L6}
        \addplot[
            color=orange,
            mark=diamond*,
            thick,
            dashed
        ] coordinates {
            (2, 5344.21)
            (4, 5789.56)
            (6, 6947.47)
            (8, 7642.21)
        };
        \addlegendentry{L7}
        \end{axis}
        \end{tikzpicture}
        \caption{Execution Time evolution with Compartments}
        \label{fig:L5_veh_time}
    \end{subfigure}
    \caption{Evolution of the Objective and Execution Time with the Number of Compartments for Instances L4 to L7}
    \label{fig:L4L7_evolution}
\end{figure}
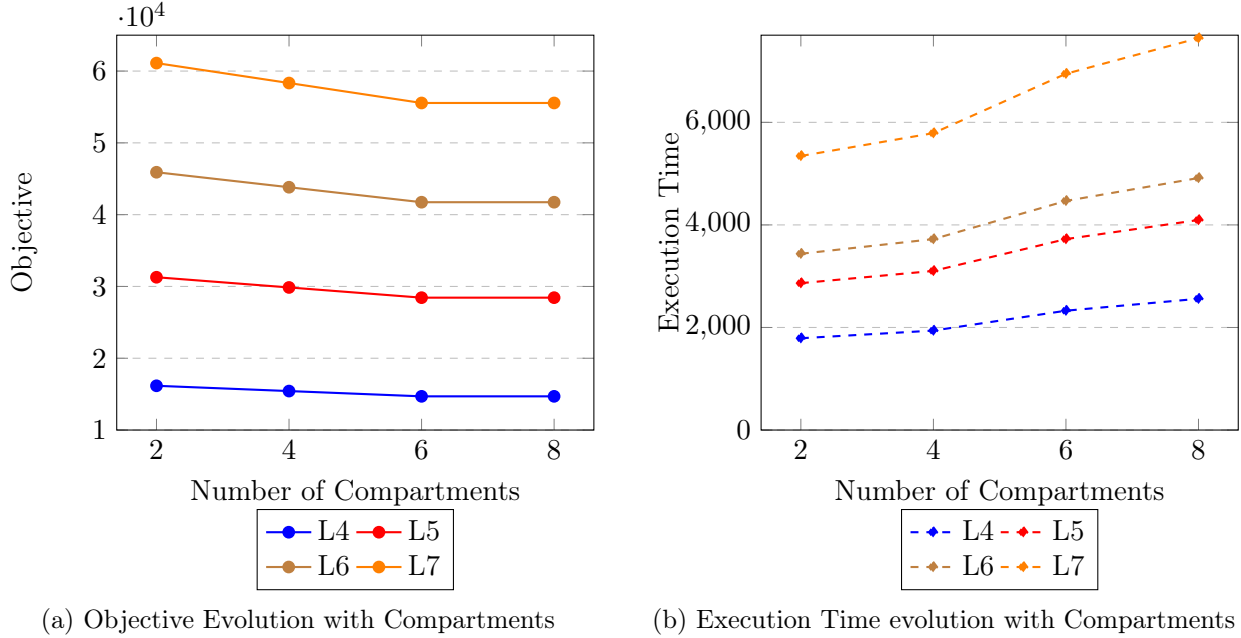

The gain in solution quality with more compartments comes at the cost of increased execution time, particularly for larger instances (see Figure \ref{fig:L4L7_evolution}). Average solve times rise from 2{,}592.00 seconds (2 compartments) to 3{,}421.44 seconds (8 compartments). This trade-off reflects the added complexity of handling finer-grained compartmental loading, which increases the model's size, the solver's burden, and memory requirements.

\subsubsection{Impact of Time Windows}

Consider the impact of the number of time windows on the total distance, number of vehicles, and execution time. Table \ref{tab:Impact of Time Windows} compares two different time window settings: Single Time Window (STW) and Multiple Time Windows (MTW) per customer.

\begin{table}[t!]
\centering
\scalebox{0.7}{
\caption{Impact of Time Windows}
\label{tab:Impact of Time Windows}
\begin{tabular}{ccc|ccc|ccc}
    \toprule
    \multirow{2}[4]{*}{\textbf{\#Inst}} & \multirow{2}[4]{*}{\textbf{\#Cust}} & \multirow{2}[4]{*}{\textbf{TL}} & \multicolumn{3}{c|}{\textbf{STW}} & \multicolumn{3}{c}{\textbf{MTW}} \\
\cmidrule{4-9}          &       &       & \textbf{Obj} & \textbf{\#Veh} & \textbf{Time} & \textbf{Obj} & \textbf{\#Veh} & \textbf{Time} \\
    \midrule
    \textbf{M1} & 10    & 21{,}292.82 & 1{,}445.53 & 1     & 2.67  & 1{,}363.71 & 1     & 3.47 \\
    \textbf{M2} & 15    & 33{,}835.52 & 1{,}688.41 & 1     & 30.16 & 1{,}592.84 & 1     & 39.21 \\
    \textbf{M3} & 20    & 47{,}658.39 & 2{,}963.12 & 2     & 49.43 & 2{,}795.40 & 2     & 64.26 \\
    \textbf{M4} & 25    & 59{,}608.06 & 3{,}205.76 & 2     & 70.42 & 3{,}024.30 & 2     & 91.54 \\
    \textbf{M5} & 30    & 71{,}263.73 & 3{,}471.67 & 2     & 128.99 & 3{,}275.16 & 2     & 167.68 \\
    \textbf{M6} & 35    & 84{,}511.10 & 4{,}802.82 & 3     & 299.60 & 4{,}530.97 & 3     & 389.47 \\
    \textbf{M7} & 40    & 94{,}936.14 & 5{,}265.38 & 3     & 428.58 & 4{,}967.34 & 3     & 557.16 \\
    \midrule
    \textbf{Avg} & \textbf{25} & \textbf{59{,}015.11} & \textbf{3{,}263.24} & \textbf{2} & \textbf{144.26} & \textbf{3{,}078.53} & \textbf{2} & \textbf{187.54} \\
    \midrule
    \textbf{L1} & 50    & 117{,}049.69 & 5{,}666.88 & 3     & 968.67 & 5{,}346.12 & 3     & 1{,}259.27 \\
    \textbf{L2} & 60    & 139{,}881.54 & 7{,}507.22 & 4     & 1{,}065.54 & 7{,}082.29 & 4     & 1{,}385.20 \\
    \textbf{L3} & 75    & 200{,}675.54 & 13{,}145.22 & 6     & 1{,}278.64 & 12{,}401.15 & 6     & 1{,}662.24 \\
    \textbf{L4} & 100   & 259{,}114.57 & 15{,}575.19 & 7     & 1{,}790.10 & 14{,}693.57 & 7     & 2{,}327.13 \\
    \textbf{L5} & 200   & 467{,}837.82 & 30{,}144.44 & 14    & 2{,}864.16 & 28{,}438.15 & 13    & 3{,}723.41 \\
    \textbf{L6} & 300   & 725{,}175.74 & 44{,}237.70 & 21    & 3{,}436.99 & 41{,}733.68 & 20    & 4{,}468.09 \\
    \textbf{L7} & 400   & 953{,}081.54 & 58{,}889.68 & 27    & 5{,}344.21 & 55{,}556.30 & 25    & 6{,}947.47 \\
    \midrule
    \textbf{Avg} & \textbf{169} & \textbf{408{,}973.78} & \textbf{25{,}023.76} & \textbf{12} & \textbf{2{,}392.61} & \textbf{23{,}607.32} & \textbf{11} & \textbf{3{,}110.40} \\
    \bottomrule
    \end{tabular}%
}
\end{table}%

The results show that allowing multiple time windows improves solution quality across both medium and large instances. The average objective value decreased from 25{,}023.76 to 23{,}607.32, and the average number of vehicles dropped from 12 to 11 when moving from single to multiple time windows. This improvement is attributed to the increased temporal flexibility, which helps the solver find more efficient routes and better customer groupings. Notably, medium instances benefit slightly more, with the cost decreasing by about 5.7\%, suggesting that even moderate problem sizes can capitalize on time flexibility.

On the other hand, the use of multiple time windows significantly increases execution times, reflecting the added complexity of the problem. The average solve time rose from 2{,}392.61 seconds under the single-window model to 3{,}110.40 seconds with MTW. This pattern holds across most instances, particularly in the larger cases where additional time windows increase the size of the solution space and lead to longer solver runtimes. Still, for applications where cost and vehicle minimization are priorities (e.g., logistics or delivery), the MTW configuration offers a favorable trade-off between solution quality and computational effort, especially when modern solvers are used with sufficient computing resources.

\subsubsection{Impact of Clustering}

\begin{table}[t!]
\centering
\scalebox{0.75}{
\caption{Impact of Clustering}
\label{Impact of clustering}%
    \begin{tabular}{c|c|c|ccc|ccc|ccc}
    \toprule
    \multirow{2}[4]{*}{\textbf{\#Inst}} & \multirow{2}[4]{*}{\textbf{\#Cust}} & \multirow{2}[4]{*}{\textbf{TL}} & \multicolumn{3}{c|}{\textbf{15-20 Clusters}} & \multicolumn{3}{c|}{\textbf{20-25 Clusters}} & \multicolumn{3}{c}{\textbf{25-30 Clusters}} \\
\cmidrule{4-12}          &       &       & \textbf{Obj} & \textbf{\#Veh} & \textbf{Time} & \textbf{Obj} & \textbf{\#Veh} & \textbf{Time} & \textbf{Obj} & \textbf{\#Veh} & \textbf{Time} \\
    \midrule
    \textbf{L1} & 50    & 117{,}049.69 & 5{,}613.42 & 4     & 1{,}177.42 & 5{,}346.12 & 3     & 1{,}259.27 & 5{,}346.12 & 3     & 1{,}422.97 \\
    \textbf{L2} & 60    & 139{,}881.54 & 7{,}436.40 & 5     & 1{,}295.16 & 7{,}082.29 & 4     & 1{,}385.20 & 7{,}082.29 & 4     & 1{,}565.27 \\
    \textbf{L3} & 75    & 200{,}675.54 & 13{,}021.20 & 7     & 1{,}554.19 & 12{,}401.15 & 6     & 1{,}662.24 & 12{,}401.15 & 6     & 1{,}878.33 \\
    \textbf{L4} & 100   & 259{,}114.57 & 15{,}428.25 & 8     & 2{,}175.87 & 14{,}693.57 & 7     & 2{,}327.13 & 14{,}693.57 & 7     & 2{,}629.66 \\
    \textbf{L5} & 200   & 467{,}837.82 & 29{,}860.06 & 14    & 3{,}481.39 & 28{,}438.15 & 13    & 3{,}723.41 & 28{,}438.15 & 13    & 4{,}207.45 \\
    \textbf{L6} & 300   & 725{,}175.74 & 43{,}820.37 & 21    & 4{,}177.66 & 41{,}733.68 & 20    & 4{,}468.09 & 41{,}733.68 & 20    & 5{,}048.94 \\
    \textbf{L7} & 400   & 953{,}081.54 & 58{,}334.12 & 26    & 6{,}495.88 & 55{,}556.30 & 25    & 6{,}947.47 & 55{,}556.30 & 25    & 7{,}850.64 \\
    \midrule
    \textbf{Avg} & \textbf{169} & \textbf{408{,}973.78} & \textbf{24{,}787.69} & \textbf{12} & \textbf{2{,}908.22} & \textbf{23{,}607.32} & \textbf{11} & \textbf{3{,}110.40} & \textbf{23{,}607.32} & \textbf{11} & \textbf{3{,}514.75} \\
    \bottomrule
    \end{tabular}%
}
\end{table}%

Recall that since because the large-scale instances (second class) cannot be solved with B\&P, RS-B\&P uses small clusters of 20 to 25 customers. Consider the impact of the clusters on the total distance, the number of vehicles, and the execution time. Table \ref{Impact of clustering} highlights the impact of the cluster size on solving the large instances. The results compare three clustering approaches: 15–20, 20–25, and 25–30 clusters. Clustering consistently reduces computational time, particularly for larger instances, while maintaining solution quality within an acceptable range. As the number of clusters increases, the objective values stabilize, indicating that finer clustering improves solution accuracy without significantly increasing the number of vehicles.

The total solution time increases slightly with an increase in the number of clusters due to a greater number of subproblems. Clustering also tends to use fewer or the same number of vehicles, enhancing resource efficiency. In particular, with 20–30 clusters, we obtain a practical tradeoff between speed and solution quality.

\subsection{Managerial Insights}\label{section5:ManagerialInsights}

This research was inspired by a real-world challenge faced by a major U.S. company, which relied heavily on an experienced operator with two decades of expertise to manually plan delivery routes. This traditional approach involved crafting routes based on intuition and experience, followed by manual checks to ensure that loading constraints for their multi-compartment vehicles were respected. However, the routes designed by the operator, while practical, tended to be shorter and less efficient compared to those obtained using the B\&P approach. By leveraging decomposition methods, the approach proposed in this paper combines and refines partial routes identified manually, enabling fuller vehicle utilization and significant reductions in total distance traveled, as well as the number of vehicles.

From a practical standpoint, the benefits of this transition are substantial. For the plant under study, the optimized solutions demonstrate a 20\% reduction in the number of vehicles required, resulting in significant operational savings. Additionally, the total distance driven is reduced by approximately 7,000 kilometers, which translates not only into fuel savings but also into lower emissions and maintenance costs. When deployed in practice, the industrial partner estimates \$250,000 in annual gains for this single facility, representing a significant operational impact. Importantly, the industrial partner has expressed strong interest in deploying the same optimization approach across its other plants, which would further amplify the impact at the network level. This underscores the scalability and practical relevance of the RS-B\&P approach, making it highly applicable for complex logistics environments.

Beyond direct cost savings, the adoption of optimization techniques is known to offer strategic, tactical, and operational advantages. Automated and systematic route planning reduces reliance on individual operator expertise, mitigating risks related to knowledge loss or subjective decision-making. It also increases responsiveness and adaptability, enabling companies to quickly adjust routes based on changing demand or operational constraints by running several what-if scenarios. Furthermore, improved vehicle utilization fosters sustainability goals by minimizing unnecessary trips and lowering carbon footprints. Overall, this research underscores the transformative power of integrating advanced optimization into logistics operations, providing companies with both immediate financial benefits and a foundation for long-term competitive advantage in an increasingly complex supply chain landscape.

\section{Conclusions}\label{section:6}

This study addressed the multi-compartment vehicle routing problem with multiple time windows (MCVRPMTW), a generalization of the classical VRPTW that incorporates complex compartment-related constraints and multiple client-specific time windows. It considered key features such as compartment size flexibility, item-to-compartment compatibility, and item-to-item compatibility, as well as the presence of multiple time windows per customer. To solve this problem, the study developed a branch-and-price algorithm enhanced with a labeling algorithm for the pricing problem. Several acceleration strategies were integrated to improve computational efficiency, including symmetry reduction in label extensions, column selection to improve dual stability, and a more effective branching process through primal diving heuristics. Additionally, the study introduced a new set of benchmark instances inspired by a real-world large-scale MCVRPMTW case provided by our industrial partner. To tackle these large-scale instances, the solution method includes a rolling space algorithm that combines our B\&P with clustering techniques. Extensive computational experiments demonstrated the effectiveness of the proposed approach and provided valuable managerial insights for practical implementations.

\section*{Acknowledgments}

This research was partly supported by the NSF AI Institute for Advances in Optimization (Award 2112533).

\begin{spacing}{1}
\typeout{}
\bibliographystyle{apalike}
\bibliography{References.bib}
\end{spacing}    

\vspace{-0.3cm}

\begin{appendices}
\renewcommand{\thesection}{\Alph{section}}
\section{Compact Formulation} 
\label{A}

\begin{table}[H]
\caption{Problem Notation}
\label{Problem Notation}
\centering
\begin{adjustbox}{width=\textwidth, center}
\begin{tabular}{ll}
\hline
Notation & Definition \\
\hline
\text{Sets} \quad &\\
$\mathcal{N}$ \quad & \text{Set of all clients nodes} $\ n \in \mathcal{N}$\\
$\mathcal{V}$ \quad & \text{Set of all nodes of the network, including the source and destination nodes in addition} ($\mathcal{V} = \mathcal{N}  \cup \{\sigma,\delta\}$) \\
$\mathcal{A}$ \quad &  \text{Set of possible connections (arcs)} \\
$\mathcal{S}$ \quad & \text{Subset of vertices in the sub-tour} \\
$\mathcal{K}$ \quad & \text{Set of tank trucks $k \in \mathcal{K}$} \\
$\mathcal{D}$ \quad & \text{Set of days} $d \in \mathcal{D}$\\
$\Theta_{n}$ \quad& \text{Set of time windows corresponding to client} $n \in \mathcal{N}$ \\
\text{Parameters} \quad &\\
$\boldsymbol{c}_{ij}$ \quad & \text{Cost of visiting node}  $j$  \text{after node}  $i$ \\
$\boldsymbol{Q}_n$ \quad & \text{Order load for client}  $n \in \mathcal{N}$ \\
$\boldsymbol{L}_{k}$ \quad & \text{Capacity of truck}  $k \in \mathcal{K}$ \\
$\overline{\boldsymbol{L}}_{mk}$ \quad & \text{Maximal capacity of compartment} $m \in \mathcal{M}_k$ \text{in truck} $k \in \mathcal{K}$\\
$\boldsymbol{N}_{mk}$ \quad & \text{Capacity on clients on compartment} $m \in \mathcal{M}_k$ \text{in truck} $k \in \mathcal{K}$\\
$\boldsymbol{\theta}_n$ \quad & \text{Number of time windows for client} $n \in \mathcal{N}$ \\
$\boldsymbol{s}_n^p$ \quad & \text{Lower bound time window} $p$ {at client} $n \in \mathcal{N}$ \\
$\boldsymbol{e}_n^p$ \quad & \text{Upper bound time window} $p$ {at client} $n \in \mathcal{N}$ \\
$\boldsymbol{u}_n$ \quad & \text{Service time at client} $n \in \mathcal{N}$ \\
$\boldsymbol{t}_{ij}$ \quad & \text{Travel time between clients} $i$ \text{and} $j$\\
$\boldsymbol{b}_{nm}^k$ \quad & \text{Binary parameter equal to 1 if the client} $n \in \mathcal{N}$ \text{is compatible with compartment} $m \in \mathcal{M}_k$ \text{of vehicle} $k \in \mathcal{K}$, 0  \text{otherwise}  \\
$\boldsymbol{f}_{ij}$ \quad & \text{Binary parameter equal to 1 if client} $i \in \mathcal{N}$ \text{and client} $j \in \mathcal{N}$ \text{are compatible}, 0  \text{otherwise} \\
$\boldsymbol{M}$ \quad & \text{Sufficiently large number} \\
$\boldsymbol{T}$ \quad & \text{Maximum allowable working time for a driver on day} $d \in \mathcal{D}$ \\
$\boldsymbol{dist}(i,j)$ \quad & \text{Travel distance between nodes} $i,j \in \mathcal{V}$ \\
$\boldsymbol{Dist}$ \quad & \text{Maximum allowable driving distance for a driver on day} $d \in \mathcal{D}$ \\
\text{Decision Variables} \quad &\\
$x^{kd}_{ij}$ \quad & \text{Binary variable equal to} 1 \text{if} \text{the tank truck}  $k \in \mathcal{K}$   \text{visits node}  $j$  \text{on day}  $d \in \mathcal{D}$  \text{after visiting the node} $i$, 0  \text{otherwise} \\
$y^{k}_{mn}$ \quad & \text{Real variable for the loaded quantity for client} $n \in \mathcal{N}$ \text{is on compartment} $m \in \mathcal{M}_k$ \text{of truck} $k \in \mathcal{K}$ \\
$z^{k}_{mn}$ \quad & \text{Binary variable equal to} 1  \text{if the product corresponding to client}  $n \in \mathcal{N}$  \text{is loaded in compartment}  $m \in \mathcal{M}_k$  \text{of tank truck}  $k \in \mathcal{K}$\\
$a_n^k$ \quad & \text{Real variable for the arrival time of vehicle} $k \in \mathcal{K}$ \text{at client} $n \in \mathcal{N}$ \\
$q_n^k$ \quad & \text{Real variable for the departure time of vehicle} $k \in \mathcal{K}$ \text{from client} $n \in \mathcal{N}$\\
$v^{p}_{n}$ \quad & \text{Binary variable equal to} 1  \text{if and only if the client}  $n \in \mathcal{N}$  \text{is served within its time window}  $p \in \Theta_n$ \\
\hline
\end{tabular}
\end{adjustbox}
\end{table}

\footnotesize{
\begingroup
\allowdisplaybreaks
\centering
\begin{align}
\label{Objective}
\min & \sum_{(i,j)\in \mathcal{A}}\sum_{k\in \mathcal{K}}\sum_{d\in \mathcal{D}} \boldsymbol{c}_{ij}x_{ij}^{kd} &&\\
\notag s.t. \hspace{2mm} & \textrm{One Trailer Truck and One Delivery Day per Customer:} &&\\
\label{One Trailer per Customer} &\sum_{(i,n)\in \mathcal{A}}\sum_{k \in \mathcal{K}}\sum_{d \in \mathcal{D}}x_{in}^{kd}=1 && \forall n \in \mathcal{N}\\
& \notag \textrm{Multi-Compartment Constraints:}\\
\label{Loading 1} & z_{mn}^k \leq \boldsymbol{b}_{mn}^k && \forall n \in \mathcal{N}, k \in \mathcal{K}, m \in \mathcal{M}_k \\
\label{Loading 2} & \sum_{m \in \mathcal{M}_k} y_{mn}^k = \boldsymbol{Q}_n \sum_{(i,n)\in \mathcal{A}}\sum_{d \in \mathcal{D}}x_{in}^{kd} && \forall n \in \mathcal{N}, k \in \mathcal{K} \\
\label{Loading 3} & \sum_{n\in \mathcal{N}}\sum_{m \in \mathcal{M}_k} y_{mn}^k \leq 
\boldsymbol{L}_k && \forall k \in \mathcal{K} \\
\label{Loading 4} & \sum_{n\in \mathcal{N}} y_{mn}^k \leq 
\overline{\boldsymbol{L}}_{mk} \sum_{n\in \mathcal{N}} z_{mn}^k && \forall k \in \mathcal{K}, m \in \mathcal{M}_k \\
\label{Loading 6} & y_{mn}^k \leq 
\boldsymbol{Q}_n z_{mn}^k && \forall n \in \mathcal{N},  k \in \mathcal{K}, m \in \mathcal{M}_k \\
\label{Loading 7} & \sum_{n \in \mathcal{N}} z_{mn}^k \leq \boldsymbol{N}_{mk} && \forall k \in \mathcal{K}, m \in \mathcal{M}_k \\
\label{Loading 8} & z_{mi}^k + z_{mj}^k \leq 1 + \boldsymbol{f}_{ij} && \forall i,j \in \mathcal{N}, k \in \mathcal{K}, m \in \mathcal{M}_k \\
& \notag \textrm{Multiple Time Windows:}\\
\label{Windows 1} &  q_n^k \geq a_n^k + \boldsymbol{u}_n - \boldsymbol{M} (1-\sum_{m \in \mathcal{M}_k}z_{mn}^k) && \forall n \in \mathcal{N}\\
\label{Windows 2} &  a_j^k \geq q_i^k + \boldsymbol{t}_{ij} - \boldsymbol{M} (1-x_{ij}^{kd}) && \forall (i,j) \in \mathcal{A}, k \in \mathcal{K} \\
\label{Windows 3} &  a_j^k \leq q_i^k + \boldsymbol{t}_{ij} + \boldsymbol{M} (1-x_{ij}^{kd}) && \forall (i,j) \in \mathcal{A}, k \in \mathcal{K} \\
\label{Windows 4} &  a_i^k \geq s_n^p - \boldsymbol{M}(1-\sum_{m \in \mathcal{M}_k}z_{mn}^k) - \boldsymbol{M}(1-v_i^p) && \forall i \in \mathcal{N}, p \in \Theta_n, k \in \mathcal{K} \\
\label{Windows 5} &  a_i^k \leq e_n^p + \boldsymbol{M}(1-\sum_{m \in \mathcal{M}_k}z_{mn}^k) + \boldsymbol{M}(1-v_i^p) && \forall i \in \mathcal{N}, p \in \Theta_n, k \in \mathcal{K} \\
\label{Windows 6} &\sum_{p \in \Theta}v_{n}^{p}=1 && \forall n \in \mathcal{N}\\   
& \notag \textrm{Flow Conservation \& Subtour Elimination:}\\
\label{flow 1} & \sum_{i \in \mathcal{V}} x^{k}_{ij} - \sum_{i \in \mathcal{V}} x^{k}_{ji} = 
     \begin{cases}
       -1 & \text{if $j=\sigma$} \\
       0 & \text{if $j\in \mathcal{V} \setminus \{\sigma,\delta\}$} \\
       1 & \text{if $j=\delta$}
      \end{cases}  
     && \forall k \in \mathcal{K}  \\
\label{flow 2} & \sum_{i,j \in \mathcal{S}} x^{k}_{ij} \leq |\mathcal{S}|-1 
     &&  \forall k \in \mathcal{K} \hspace{1mm} (\mathcal{S} \subset \mathcal{V}, 2 \leq |\mathcal{S}| \leq |\mathcal{V}|-2)\\
& \notag \textrm{Driver Regulations:} && \\
\label{Regulations 1} & \sum_{(i,j) \in \mathcal{A}} \boldsymbol{t}_{ij} x_{ij}^{kd} \leq \boldsymbol{T} && \forall k \in \mathcal{K}, d \in \mathcal{D}\\
\label{Regulations 2} & \sum_{d \in \mathcal{D}}\sum_{(i,j) \in \mathcal{A}} \boldsymbol{t}_{ij} x_{ij}^{kd} \leq \boldsymbol{T} && \forall k \in \mathcal{K}\\
\label{Regulations 3} & \sum_{(i,j) \in \mathcal{A}} \boldsymbol{dist}_{ij} x_{ij}^{kd} \leq \boldsymbol{Dist} && \forall k \in \mathcal{K}, d \in \mathcal{D}\\
& \notag \textrm{Conditions:}\\
\label{Conditions 1} & x_{ij}^{kd}\in\{0,1\} && \forall (i,j) \in \mathcal{A}, k \in \mathcal{K}, d \in \mathcal{D} \\
\label{Conditions 2} & y^{k}_{mn}\geq 0 && \forall n \in \mathcal{N}, k\in \mathcal{K}, m \in \mathcal{M}_k \\
\label{Conditions 3} & z^{k}_{mn}\in \{0,1\} && \forall n \in \mathcal{N}, k\in \mathcal{K}, m \in \mathcal{M}_k\\
\label{Conditions 4} & a^{k}_{n}\geq 0 && \forall n \in \mathcal{N}, k\in \mathcal{K}\\
\label{Conditions 5} & q^{k}_{n}\geq 0 && \forall n \in \mathcal{N}, k\in \mathcal{K}\\
\label{Conditions 7} & v^{p}_{n}\in \{0,1\} && \forall n \in \mathcal{N}, p \in \Theta
\end{align}
\endgroup
}

\vspace{-0.3cm}

\normalsize The objective function (\ref{Objective}) minimizes the total cost. Constraints (\ref{One Trailer per Customer}) ensure each customer is served by a single truck. Constraints (\ref{Loading 1}) check whether a client can be loaded on a given compartment on a given truck. Constraints (\ref{Loading 2}) ensure the load corresponding to a given client on the truck serving this customer is equal to the load ordered by the customer. Constraints (\ref{Loading 3}) control the truck capacity. Constraints (\ref{Loading 4}) control the compartment capacities. Constraints (\ref{Loading 6}) ensure that a client is not loaded on a given truck compartment unless allowed. Constraints (\ref{Loading 7}) ensure that each compartment serves at most $\boldsymbol{N}_{mk}$. Constraints (\ref{Loading 8}) control the clients' compatibility. Constraints (\ref{Windows 1}) ensure that the departure from client $n$ is at least equal to the arrival time at client $n$, plus the waiting time and the service time at customer $n$, only if client $n$ is assigned to vehicle $k$. Constraints (\ref{Windows 2}) and (\ref{Windows 3}) mean that the arrival time at client $j$ is equal to the departure time from client $i$, plus the travel time $\boldsymbol{t}_{ij}$ only if the arc $(i,j)$ is assigned to vehicle $k$. Constraints (\ref{Windows 4}) and (\ref{Windows 5}) imply that the arrival time plus the waiting time of vehicle $k$ at client $i$ is within the time window $p$ only if the client is assigned to vehicle $k$ and time window $p$ is chosen. Constraints (\ref{Windows 6}) mean that exactly one single time window is chosen for each customer $n$. Constraints (\ref{flow 1}) and (\ref{flow 2}) are the flow conservation and subtour elimination constraints. Constraints (\ref{Regulations 1}) ensure that the drivers do not exceed a maximum working time each day. Constraints (\ref{Regulations 2}) ensure that the drivers do not exceed a maximum working time during the time horizon. Constraints (\ref{Regulations 3}) ensure that the drivers do not exceed a maximum driving distance each day. Constraints (\ref{Conditions 1}) to (\ref{Conditions 7}) are the conditions on the problem variables.

\section{Attribute Settings}
\label{C}
Table \ref{Summary of Attribute Settings} summarizes the features of the tested cases. For instance, the fifth line indicates the item-to-compartment compatibility for the 6-compartment case. Items 1 and 2 are compatible with all compartments, items  3 to 6 are compatible with the first 4 compartments (side compartments), while items 7 to 8 are compatible with compartments 5 and 6 (back compartments). Similarly, the sixth line indicates the item-to-item compatibility for the same case. Items 1 and 2 are compatible with all other items. Items with attributes 3, 5, and 7 are compatible with all other odd attributes. Items 4, 6, and 8 with attributes are compatible with all other even attributes.

\begin{table}[H]
\centering
\caption{Summary of Attribute Settings}
\label{Summary of Attribute Settings}
\scalebox{0.8}{
    \begin{tabular}{c|c|c|c|c|c|c|c|c|c}
    \toprule
    \textbf{Case} & \textbf{Attribute} & \textbf{1} & \textbf{2} & \textbf{3} & \textbf{4} & \textbf{5} & \textbf{6} & \textbf{7} & \textbf{8} \\
    \midrule
    \multirow{2}[4]{*}{\textbf{2}} & \textbf{i-to-c} & All   & All   & \{1\}     & \{1\}     & \{1\}     & \{2\}     & \{2\}     & \{2\} \\
\cmidrule{2-10}          & \textbf{i-to-i} & All   & All   & \{1,5,7\} & \{2,6,8\} & \{1,3,7\} & \{2,4,8\} & \{1,3,5\} & \{2,4,6\} \\
    \midrule
    \multirow{2}[4]{*}{\textbf{4}} & \textbf{i-to-c} & All   & All   & \{1,2\} & \{1,2\} & \{1,2\} & \{3,4\} & \{3,4\} & \{3,4\} \\
\cmidrule{2-10}          & \textbf{i-to-i} & All   & All   & \{1,5,7\} & \{2,6,8\} & \{1,3,7\} & \{2,4,8\} & \{1,3,5\} & \{2,4,6\} \\
    \midrule
    \multirow{2}[4]{*}{\textbf{6}} & \textbf{i-to-c} & All   & \{1,2\}   & \{1,2,3,4\} & \{1,2,3,4\} & \{1,2,3,4\} & \{1,2,3,4\} & \{3,4,5,6\} & \{3,4,5,6\} \\
\cmidrule{2-10}          & \textbf{i-to-i} & All   & All  & \{1,5,7\} & \{2,6,8\} & \{1,3,7\} & \{2,4,8\} & \{1,3,5\} & \{2,4,6\} \\
    \midrule
    \multirow{2}[4]{*}{\textbf{8}} & \textbf{i-to-c} & All   & \{1,2\}   & \{1,2,3,4,5,6\} & \{1,2,3,4,5,6\} & \{1,2,3,4,5,6\} & \{1,2,3,4,5,6\} & \{5,6,7,8\} & \{5,6,7,8\} \\
\cmidrule{2-10}          & \textbf{i-to-i} & All   & All  & \{1,5,7\} & \{2,6,8\} & \{1,3,7\} & \{2,4,8\} & \{1,3,5\} & \{2,4,6\} \\
    \bottomrule
    \end{tabular}%
}
\end{table}

\end{appendices}

\end{document}